\documentclass[reqno,12pt]{amsart}
\usepackage{amsmath, amssymb, amsthm, amsfonts, mathtools} 
\usepackage[english]{babel}
\usepackage{bbm}
\usepackage{graphicx}
\usepackage{url}
\usepackage{soul}
\usepackage{epstopdf}
\usepackage[ruled,vlined]{algorithm2e}
\usepackage[a4paper,bindingoffset=0.5cm,left=2cm,right=2cm,top=2.5cm,bottom=2cm,footskip=.8cm]{geometry}
\usepackage{rotating}
\usepackage{amsbsy,enumerate}
\usepackage{comment}
\usepackage{mathrsfs} 
\usepackage{enumitem}
\usepackage{cancel}
\usepackage{subcaption}

\newcommand{\E}{\mathbb{E}}

\newcommand{\bs}{\boldsymbol}

\newcommand{\vb}{\vspace{3.2mm}}
\renewcommand{\hat}{\widehat}

\newcommand{\vertiii}[1]{{\left\vert\kern-0.25ex\left\vert\kern-0.25ex\left\vert #1 \right\vert\kern-0.25ex\right\vert\kern-0.25ex\right\vert}}
\allowdisplaybreaks
\usepackage{hyperref}

\newcommand{\Var}{\mathbb{V}{\rm ar}}
\newcommand{\Cov}{\mathbb{C}{\rm ov}}

\newcommand{\asarrow}{\:\to_{\rm as}\:}

\allowdisplaybreaks

\newtheorem{lemma}{Lemma}

\newtheorem{assumption}{Assumption}
\newtheorem{theorem}{Theorem}
\newtheorem{remark}{Remark}

\newtheorem{proposition}{Proposition}

\usepackage[utf8]{inputenc}
\usepackage{type1cm}         
\usepackage{multicol}        % used for the two-column index
\usepackage[bottom]{footmisc}% places footnotes at page bottom

\usepackage{newtxtext}       % 
\usepackage{newtxmath}       % selects Times Roman as basic font

\usepackage{pgfplots}

\usepackage{tikz}
\usetikzlibrary{automata,positioning}
\usetikzlibrary{decorations.pathreplacing,decorations.markings,shapes.geometric}
\usetikzlibrary{shapes,arrows}
\usetikzlibrary{backgrounds,calc,positioning}

\usepackage{pgfplots}
\pgfplotsset{compat=1.9}

\usetikzlibrary{chains,shapes.multipart}
\usetikzlibrary{shapes,calc}
\usetikzlibrary{automata,positioning}

\usepackage{tikz}
\usetikzlibrary{automata,positioning}
\usetikzlibrary{decorations.pathreplacing,decorations.markings,shapes.geometric}
\usetikzlibrary{shapes,arrows}
\usetikzlibrary{backgrounds,calc,positioning}

\newcommand{\michel}[1]{{\color{black}#1}}
\newcommand{\hritika}[1]{{\color{black}#1}}
\newcommand{\liron}[1]{{\color{black}#1}}

\begin{document}

	\title[Uncovering the topology of a queueing network from population data]{Uncovering the topology of an \\infinite-server queueing network\\ from population data}
\author[{H. Gupta, M. Mandjes, L. Ravner {\tiny and} J. Wang}]{Hritika Gupta, Michel Mandjes, Liron Ravner {\tiny and} Jiesen Wang}

%\textcolor{blue}{Alternative title suggestions:\\ Uncovering the structure of a network of infinite server queues from node count data;\\
%Uncovering the topology of a network of infinite server queues from population data; \\ Parametric inference for infinite server queueing networks with partial data}

	\begin{abstract}
			This paper studies statistical inference in a network of infinite-server queues, with the aim of estimating the underlying parameters (routing matrix, arrival rates, parameters pertaining to the service times) using observations of the network population vector at Poisson time points. We propose a method-of-moments estimator and establish its consistency. The method relies on deriving the covariance structure of different nodes at different sampling epochs. Numerical experiments demonstrate that the method yields accurate estimates, even in settings with a large number of parameters. Two model variants are considered: one that assumes a known parametric form for the service-time distributions, and a model-free version that does not require such assumptions.

\vb


\vb

\noindent
{\sc Affiliations.} 
HG is with School of Mathematics and Statistics, The University of Melbourne, Victoria 3010, Australia. MM is  with Mathematical Institute, Leiden University, Leiden, 
The Netherlands; he is also affiliated with (a)~Korteweg-de Vries Institute for Mathematics, University of Amsterdam, Amsterdam, The Netherlands, (b)~E{\sc urandom}, Eindhoven University of Technology, Eindhoven, The Netherlands, (c)~Amsterdam Business School, Faculty of Economics and Business, University of Amsterdam, Amsterdam, The Netherlands.  
LR is with Department of Statistics, University of Haifa,
Israel.
JW is with Korteweg-de Vries Institute for Mathematics, University of Amsterdam, Amsterdam, The Netherlands.

\noindent Date: {\it \today}.

\vb

\noindent
{\sc Acknowledgments.} 
\liron{The authors wish to thank the associate editor and the reviewer for their insightful and highly useful feedback, which helped improve this paper significantly.}
MM's and JW's research has been funded by the European Union’s Horizon 2020 research and innovation programme under the Marie Sklodowska-Curie grant agreement no.\ 945045, and by the NWO Gravitation project {\tiny NETWORKS} under grant agreement no.\ 024.002.003. LR's research was supported by the Israel Science Foundation (ISF), grant no.\ 1361/23. \includegraphics[height=1em]{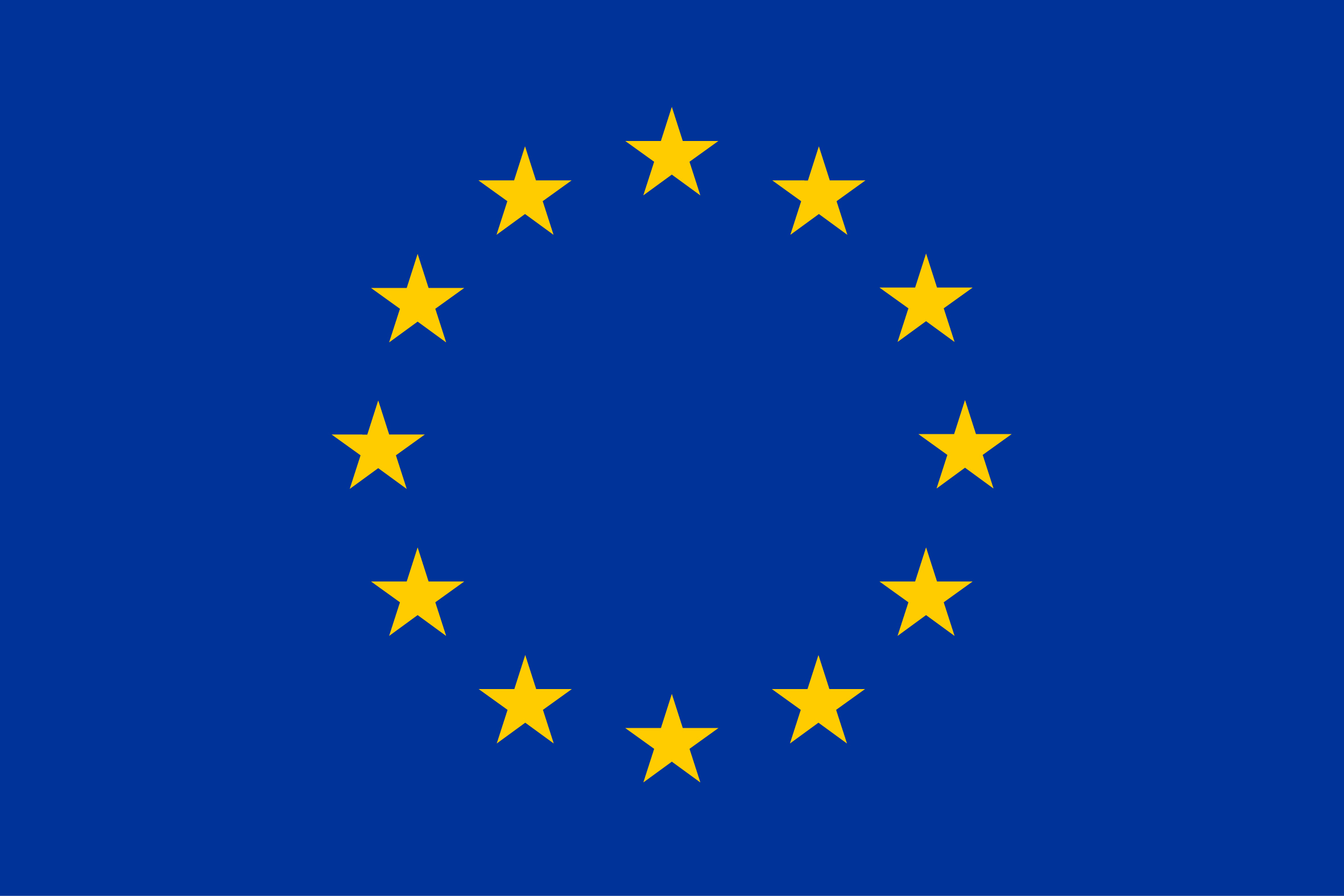} HG's research has been funded by the Melbourne Research Scholarship, and by the Australian Research Council Industrial Transformation Training Centre in Optimisation Technologies, Integrated Methodologies, and Applications (OPTIMA), Project ID IC200100009.

\vb

\noindent
{\sc Email.} \url{hritikag@student.unimelb.edu.au}, \:\url{m.r.h.mandjes@math.leidenuniv.nl}, 

\noindent \url{lravner@stat.haifa.ac.il}, \: \url{j.wang2@uva.nl}

	\end{abstract}

 \maketitle

%{\small Things to think about:
%\begin{itemize}
%\item[--] Step (iii) in the proof of Proposition 2: I don't understand (anymore) how uniqueness follows. \textcolor{blue}{LR: Added some thoughts about this is Section 4.2.} 
% \item[--] In \S 4.4 LR wrote: `Perhaps we need to be more careful with the Lipschitz constants here? Do we need the classes to have the same uniform constant? I think this is true for the sum arguments, but not for the product.' MM: I thought this works out, as it is always about finitely many terms, so one can always work with the largest of the Lipschitz constants. Perhaps I don't understand the problem. \textcolor{blue}{LR: After reading again I see your point. The argument is indeed sufficient.} 
%\item[--] MM: in the end I have gone in \S 4.6 for first explaining an approach in which the service times are dealt with in a parametric fashion, and then for an approach that works with ${\mathbb E}[G_i]$, $\mathscr{G}_i(\beta)$, and $\dot{\mathscr{G}}_i(\beta)$. I think it makes sense to present both.  \textcolor{blue}{LR: I like this approach.} 
%\end{itemize}}

 \newpage

 \section{Introduction}

In operations research, when analyzing a stochastic system, it is typically assumed that the underlying model is known and expressed in terms of one or more stochastic processes, each characterized by a set of parameters. The goal is then to evaluate specific performance metrics, which inform system design; for example, being able to compute these metrics allows for optimizing an objective function subject to constraints. However, in practice, model parameters are not known a priori and must be inferred from data. The simplest scenario occurs when the stochastic processes are directly observed, making inference a conventional statistical task. In many cases, though, these processes are not observed directly (see \cite{SJ2011}), complicating the inference process significantly. 

A prototypical example of inference based on indirect observations, in the context of operations research, is as follows. Suppose one observes the workload of a storage system at specific points in time (for instance equidistantly, or at Poisson times). Given these workload observations, can the parameters of the system's input process be estimated? This type of problems is often referred to as {\it inverse problems}: observing an endogenous
process, we wish to  estimate parameters of the corresponding exogenous process. 
Even if the observed workload process is assumed to be in stationarity, this class of inference problems is highly challenging due to the intricate dependence structure of the observed workload values. Concretely, only in very special cases it is possible to write down the likelihood pertaining to the observed data. A survey on parameter inference for queueing systems can be found in \cite{Azam}.

In this paper, we consider parameter inference in a network of {$n$} infinite-server queues. In our setup customers arrive at the network's nodes, undergo service, are routed according to a routing matrix, and eventually leave the network. Importantly, all nodes are of the infinite-server type \cite[\S 10.3]{whitt}, meaning that customers can be served in parallel, i.e., they do not experience any waiting time. We throughout assume  that the service-time distributions are either known or follow a specified parametric form. The objective of this paper is to develop an algorithm for uncovering the network structure using observations of the queue length vector only — specifically, to estimate the model's arrival rates and routing matrix. A secondary goal is to estimate the parameters pertaining to the service times. Finally, we aim to address the setting of noisy observations, in which each customer present at an observation time is detected with a certain (unknown) probability.

 \medskip

 The main contribution of this work concerns a procedure that is capable of estimating all model parameters, including the `observation probabilities'.  In our setup, the data available is the {$n$-dimensional} vector of queue lengths observed at Poisson times (with the sampling rate being known). Our approach is based on the {\it method of moments}. Under the mild assumption that we observe the network  in stationarity, we utilize explicit expressions for the expected values of the queue lengths and the covariances between any pair of queue lengths at two subsequent observation times (i.e.,  the inter-observation times have an exponential distribution).  This procedure results in $n^2+n$ moment equations: $n$ for the mean queue-lengths at each station and $n^2$ for the cross-moments corresponding to all pairs of queues observed at two consecutive observations.

In case the observation probabilities and the service-time distributions are known, then 
the number of moment equations matches the number of unknown parameters (namely: the $n^2$ entries of the routing matrix and the $n$ external arrival rates).
 In case each of the service-time distributions is characterized by a single parameter, we have to identify $n$ additional parameters. In this case, and in the case that the observation probabilities are not known, the above set of moment equations alone would be insufficient. However, this issue can be resolved by also including the covariances between any pair of queue-lengths observed at times {\it two} steps apart (i.e., in which the inter-observation times have an Erlang-2 distribution). We in addition present a `model-free approach', i.e., an approach which has the advantage that one does not need to know the parametric class of the service-time distributions.

 Our main result states that the resulting estimator is consistent: we establish that the estimator \liron{converges with probability one} to the true value of the parameter as the number of observations grows large. The proofs use tools presented in \cite{V2000}; most notably, it is shown that the moments (as being used in the moment equations) are estimated in a consistent manner, and that this consistency is inherited by our parameter estimator. Remarkably, the method can also recover the {\it directed} structure of distinct networks that share the same stationary distribution. This is possible because, for each pair $(i_1, i_2)$ of stations, it utilizes cross-moments between observation times, which are not necessarily symmetric in $i_1$ and $i_2$.

 We demonstrate the efficacy of our procedure through a series of experiments, which show that even for models with a substantial number of parameters, our method accurately estimates them. Specifically, the experiments highlight the ability of our estimator to recover the structure of the underlying routing matrix. %Additionally, we discuss several computational strategies that are effective in solving the moment equations efficiently.

 \medskip

We continue with a brief overview of the relevant literature, focusing initially on estimation in the context of {\it single} infinite-server queues. A foundational contribution is by Brown \cite{Brown}, who introduced a nonparametric estimator for the service-time distribution in the M/G/$\infty$ queue, relying solely on observations of arrival and departure times (i.e., without having any information about the way these arrival and departure times are coupled). This approach results in a consistent estimator, with further extensions being discussed in \cite{BNW}. Subsequently, \cite{G2} analyzed the estimation procedure from \cite{Brown}, presenting exact non-asymptotic expressions for the mean squared error as a key result. Another relevant contribution is \cite{BH}, which explores both parametric and nonparametric estimation techniques for M/G/$\infty$ systems, based on either the queue-length process or sequences of busy and idle periods; from a practical perspective, a challenge with the latter approach is that, in infinite-server queues with substantial load, busy periods can be disproportionately long, rendering the estimator inefficient. Finally, \cite{G1} investigates a similar estimation framework, assessing the maximal risk of the nonparametric estimator for the service-time distribution.

It is worth noting that there is a distinct body of literature focused on {\it single-server} queues. These systems are inherently more complex than infinite-server queues because customers {\it do} interfere with each other: unlike in infinite-server systems, waiting is a factor. Notable contributions in this area include \cite{HP}, which develops a non-parametric estimator for the service time distribution in the M/G/1 queue, and \cite{RBM}, which provides an estimator for the Laplace exponent of the input process in a Lévy-driven storage system. In the context of a network of {\it finite-server} queues \cite{SJ2011} suggest a Bayesian inference method for network in which only some of the jobs are observed. This is related to one of our model variants, with the important difference that in the framework of \cite{SJ2011} the full path of a sample job is observed, whereas in our framework only periodic observations of the network population vector are available. 

The literature on networks of infinite-server queues is relatively limited. Wichelhaus and Langrock \cite{WICH2} examine a class of feedforward networks (with a strong focus on systems of tandem type), where both the external arrival and departure processes are observed. Their objective is to infer the service-time distributions at the individual stations and the routing probabilities. This is achieved by employing cross-spectral techniques for multivariate point processes, building on the framework developed in~\cite{BRIL}. Additionally, Schweer and Wichelhaus~\cite{WICH} extend existing methods based on covariance functions, with the key result being a functional central limit theorem for the proposed estimators.

There is a compelling connection between the methodology developed in this paper and the concept of {\it causality} \cite{Pearl}. Causality refers to the relationship between causes and their effects, where a cause can potentially lead to the occurrence of an effect. This concept is fundamental across a wide range of scientific disciplines, including economics, physics, and sociology, where identifying causal relationships helps explain phenomena and predict future outcomes. Recognizing that causes must precede their effects temporally, time series data plays a crucial role in gaining insight into causality, particularly in determining the direction of causality. In the context of the technique presented in this paper, we are able to detect the model's routing matrix, which can be interpreted as a reflection of an underlying causal structure.

Another area of research within network-related inverse problems is {\it network tomography}. The primary objective, motivated by applications in large-scale communication networks, is to infer a network's internal characteristics using only end-point data. For a more detailed description of this field, we refer to \cite{RUI}.

\medskip

This paper is organized as follows. Section \ref{sec:model} defines our model and states the paper's objectives. Our approach is based on the method of moments; expressions for the required moments, in terms of the model parameters, are given in Section \ref{sec:mom}. Section \ref{sec:method} presents our estimator and the corresponding consistency proof for the case that the routing matrix and the external arrival rates are estimated (for given service-time distributions and observation probabilities, that is).
Section~\ref{sec:G_unknown} discusses an extension in which also service-time parameters are estimated; this section also includes the model-free approach, which does not require knowledge of a parametric class of the service-time distributions. In Section \ref{sec:obsprob} we  point out how one can also include the estimation of the observation probabilities. Numerical experiments are given in Section~\ref{sec:exp}.

 \section{Model and objectives}\label{sec:model}
In this paper we consider a network that consists of $n\in{\mathbb N}$ stations. Each of these stations is to be interpreted as an {\it infinite-server queue}, i.e., a service facility at which all customers present are served in parallel. Recall that infinite-server queues are not queues in the strict sense, in that customers do not interfere with each other. The network is characterized by three elements: the system's external arrivals, the customers' per-station service times, and the routing mechanism. 

\begin{itemize}
\item[$\circ$] {\it External arrivals.} The external arrival rate at station $i\in\{1,\ldots,n\}$ is $\lambda_i\geqslant 0$. The corresponding arrival processes are assumed to be Poisson processes.
    \item[$\circ$] {\it Service times.}
The per-customer service time in station $i\in\{1,\ldots,n\}$ is distributed as the non-negative, finite-mean random variable $G_i$ that is characterized via its Laplace-Stieltjes transform ${\mathscr G}_i(\cdot)$.
%and the probability density function by $g_i(\cdot)$. 
In our setup all service times are assumed independent, and the per-station service times are in addition assumed  to be identically distributed. 

\item[$\circ$] {\it Routing mechanism.}
We denote by $q_{ij}$ the probability that a customer goes to station $j$ after having been served at station $i$, where $i\in\{1,\ldots,n\}$ and $j\in\{0,1,\ldots,n\}$; here, the case $j=0$ corresponds to the customer leaving the network, so that $\sum_{j=0}^n q_{ij}=1$ for any $i\in\{1,\ldots,n\}$. We denote by $Q$ the $n\times n$ matrix whose $(i,j)$-th entry is given by $q_{ij}$.

\end{itemize}

\begin{remark}\label{rem:R1}\em
    In the setting where we infer the vector $\bs\lambda$ and the matrix $Q$, covered by Section \ref{sec:method}, we allow self-transitions: for all $i\in\{1,\ldots,n\}$, we may have that $q_{ii}>0$. In the setting that we also infer the parameters corresponding to the service-time distributions, covered by Section~\ref{sec:G_unknown}, we put $q_{ii}=0$ to avoid issues concerning identifiability. To see the necessity of this choice, observe that if $q_{ii}$ were positive and $G_i$ would have an exponential distribution with parameter $\mu_i$, then the time the customer spends at station $i$ is exponentially distributed with parameter $\mu_i(1-q_{ii})$, entailing that one can just estimate $\mu_i(1-q_{ii})$, and not $\mu_i$ and $q_{ii}$ individually. \hfill $\Diamond$
\end{remark}

 In our setup we impose the following two conditions on the model parameters. (i)~In the first place we want each queue $j$ to be `non-trivial', in that customers flow in at any station. This is guaranteed if \begin{equation}\label{eq:C1}
     \lambda_j>0, \:\: \mbox{or there is a path $(i_0, \ldots,i_k)$ such that $\lambda_{i_0}q_{i_0,i_1}\cdots q_{i_k j}>0$.}
 \end{equation} (ii)~In the second place we want each queue $j$ to be `non-explosive', in the sense that customers eventually flow out. This is guaranteed if 
 \begin{equation}\label{eq:C2}
  q_{j0}>0,\:\: \mbox{or there is a path $(i_0, \ldots,i_k)$ such that $q_{ji_0}\cdots q_{i_{k-1} i_k}q_{i_k0}>0$.}   
 \end{equation}

Throughout this paper, we consider this network of infinite-server queues to be (strongly) stationary. This entails that the {\it effective arrival rate} vector ${\boldsymbol \lambda}^{\rm eff}$ follows as the unique non-negative solution of  the {\it traffic equations}. Indeed, the effective arrival rate at station $i$ satisfies
\[\lambda_i^{\rm eff} = \lambda_i + \sum_{j=1}^n q_{ji} \,\lambda_j^{\rm eff},\]
or, in evident vector/matrix notation,
 \begin{equation}\label{eq:leff}{\boldsymbol\lambda}^{\rm eff} = (I-Q^\top)^{-1}{\boldsymbol \lambda};\end{equation}
 noting that $I-Q$ is a weakly chained diagonally dominant matrix
 (by virtue of the assumptions we imposed on the probabilities $q_{ij}$), it follows that there is a unique solution.

\medskip

In this paper we denote by ${\bs M}(t)$ the vector of queue lengths at the network's stations at time $t\geqslant 0$; indeed, its $i$-th entry corresponds to the number of customers at station $i$ at time $t$, with $i\in\{1,\ldots,n\}$. It is a classical result that the stationary distribution of the network's population size vector has Poisson marginals, where in addition the number of customers at the various stations are independent in stationarity. As a consequence of this fact, using that we start in stationarity, we have for any $t\geqslant 0$ and ${\boldsymbol m}\in {\mathbb N}_0^n$ that, for {\it loads} $\varrho_i$,
\[\mathbb{P}({\boldsymbol M}(t) = {\boldsymbol m}) = \prod_{i = 1}^n \frac{e^{-\varrho_i} \,(\varrho_i)^{m_i}}{m_i !};\]
here the per-queue loads are given by, for $i\in\{1,\ldots,n\}$, 
\begin{equation}\label{eq:loads}\varrho_i = \lambda_i^{\rm eff} \,\E[G_i].\end{equation}

The main goal of this paper is to infer the model parameters by observing the network population vector ${\bs M}(\cdot)$ at the times of a Poisson process with known parameter $\beta$. The method that we propose is based on the method of moments, motivating why in the next section we analyze moments and cross-moments pertaining to ${\bs M}(\cdot)$ at Poisson times.

\section{Moments pertaining to the network population}\label{sec:mom}
The main goal of this section is to derive expressions for various expectations that involve the network population vector. It is clear that, for any $t\geqslant 0$ and $i\in\{1,\ldots,n\}$,
\begin{equation}\label{eq:loads2}{\mathbb E}[M_i(t)] = \varrho_i,\end{equation}
by virtue of the multivariate Poisson stationary distribution in combination with the assumed stationarity. In this section, the object of interest is 
\[\E[M_j(t)M_i(t+T_\beta)]=\E[M_j(0)M_i(T_\beta)],\]
with $T_\beta$ an exponentially distributed random variable with mean $\beta^{-1}$, sampled independently of everything else. 
%The second subsection considers the situation with censored observations, i.e., each customers present is observed with a given probability.  

\subsection{Auxiliary results}
We start our computation of $\E[M_j(0)M_i(T_\beta)]$ by analyzing $P_{ij}(\beta)$, defined as the probability that a customer who {\it just entered} station $i$ resides at station $j$ after an independently sampled time interval $T_\beta$.
By an elementary conditioning argument, distinguishing the scenario that `the exponential clock $T_\beta$ rings' after the end of the customer's service time (at station $i$) and the scenario where these events happen in the opposite order, and in addition appealing to the memoryless property of the exponential distribution, we find
\begin{align*}
    P_{jj}(\beta) &= \int_0^\infty \int_0^t{\mathbb P(G_j\in{\rm d}x)} \,\beta\, e^{-\beta t}\, \sum_{k=1}^n q_{jk}\, P_{kj}(\beta)\, {\rm d}t + \int_0^\infty\int_t^\infty {\mathbb P(G_j\in{\rm d}x)}\, \beta \,e^{-\beta t}\,{\rm d}t\\
    &= \michel{\int_0^\infty {\mathbb P}(G_j\leqslant t) \,\beta\, e^{-\beta t}\,{\rm d}t \, \sum_{k=1}^n q_{jk}\, P_{kj}(\beta) + \int_0^\infty {\mathbb P}(G_j>t)\, \beta \,e^{-\beta t}\,{\rm d}t}\\
    &= \mathscr{G}_j(\beta) \, \sum_{k=1}^n q_{jk}\, P_{kj}(\beta) + (1-\mathscr{G}_j(\beta)),
\end{align*}
where in the last equation an elementary integration-by-parts argument has been used.
Similarly, for $i\not=j$ we find
\[ P_{ij}(\beta)= \mathscr{G}_i(\beta) \, \sum_{k=1}^n q_{ik}\, P_{kj}(\beta).\]
We thus have found, for any $j\in\{1,\ldots,n\}$, a system of $n$ equations in equally many unknowns. 
These relations can be compactly summarized in convenient vector/matrix notation, as follows. Define by ${\boldsymbol P}_j(\beta)$ the vector $(P_{1j}(\beta),\ldots,P_{nj}(\beta))^{\top}$. 
Then, with ${\rm diag}\{{\boldsymbol v}\}$ the diagonal matrix with the (column) vector ${\boldsymbol v}$ on the diagonal, and with ${\boldsymbol e}_j$ denoting the $j$-th unit vector, 
\begin{equation}
    \label{E1}
(I-{\rm diag}\{{\boldsymbol{\mathscr G}}(\beta)\}\,Q)\,{\boldsymbol P}_j(\beta) =  \big(I- {\rm diag}\{{\boldsymbol{\mathscr G}}(\beta)\}\big)\,{\boldsymbol e}_j. \end{equation}
Applying strict diagonal dominance, the matrix $(I-{\rm diag}\{{\boldsymbol{\mathscr G}}(\beta)\}\,Q)$ is non-singular, so that uniqueness of the solution follows (for any positive $\beta$). 
Note that in addition there is the option of the customer having left the network at the exponentially distributed time $T_\beta$. With $P_{i0}(\beta)$ denoting the corresponding probability, we evidently have
\[ P_{i0}(\beta) = 1- \sum_{j =1}^n P_{ij}(\beta).\]
The following lemma summarizes our findings so far. We define by ${\boldsymbol P}(\beta)$ the $n\times n$ matrix whose $(i,j)$-th entry is $P_{ij}(\beta)$, with $i,j\in\{1,\ldots,n\}$, and by ${\bs 1}$ the $n$-dimensional all ones vector. 
The above linear equations can be rewritten as, for $j\in\{1,\ldots,n\}$ and any $\beta>0$,
    \[{\boldsymbol P}_j(\beta) =  (I-{\rm diag}\{{\boldsymbol{\mathscr G}}(\beta)\}\,Q)^{-1}\big(I- {\rm diag}\{{\boldsymbol{\mathscr G}}(\beta)\}\big)\,{\boldsymbol e}_j,\]
    which further simplifies to the following statement. 
\begin{lemma}\label{L1}
    For any $\beta>0$,
    \begin{align}
        {\boldsymbol P}(\beta) =  (I-{\rm diag}\{{\boldsymbol{\mathscr G}}(\beta)\}\,Q)^{-1}\big(I- {\rm diag}\{{\boldsymbol{\mathscr G}}(\beta)\}\big).
    \end{align}
    In addition,
    \[{\boldsymbol P}_{0}(\beta) = {\bs 1}- {\boldsymbol P}(\beta) {\bs 1}. \]
\end{lemma}

In our computations so far, we have considered the situation that the customer {\it has just entered} station~$i$. In stationarity, however, a customer present at station $i$ has a {\it residual service time}, i.e., a service time that is sampled from the so-called {\it excess life distribution} \cite[\S V.1]{ASM}. Concretely, the excess life distribution pertaining to a customer residing at station $i$ is characterized through the cumulative distribution function 
\[ \michel{{\mathbb P}(G^{\rm (res)}_i\leqslant t)} := \frac{1}{\E[G_i]}\int_0^t \michel{{\mathbb P}(G_i>s)}\,{\rm d}s.\]
We proceed by calculating $P_{ij}^{\rm (res)}(\beta)$, being the counterpart
of $P_{ij}(\beta)$ but then starting with such a {residual} service time in station $i$ (rather than a `fresh' service time, that is). Using the same reasoning as above, we find, for $i \neq j$,
\begin{align*}P_{ij}^{\rm (res)}(\beta) &= \mathscr{G}^{\rm (res)}_i(\beta) \, \sum_{k=1}^n q_{ik}\, P_{kj}(\beta),\\
P_{jj}^{\rm (res)}(\beta) &= \mathscr{G}_j^{\rm (res)}(\beta) \, \sum_{k=1}^n q_{jk}\, P_{kj}(\beta) + (1-\mathscr{G}_j^{\rm (res)}(\beta)),\\
P_{i0}^{\rm (res)}(\beta)&=  1- \sum_{j =1}^n P^{\rm (res)}_{ij}(\beta),\end{align*}
where, using the above expression for the cumulative distribution function $G^{\rm (res)}_i(t)$, the Laplace-Stieltjes transform of the residual service time can be expressed in terms of its `fresh counterpart' via
\[\mathscr{G}_i^{\rm (res)}(\beta) = \frac{1-\mathscr{G}_i(\beta)}{\beta\, \E[G_i]}.\]
The next lemma presents these findings in a compact manner, where ${\boldsymbol P}^{\rm (res)}_j(\beta)$ and ${\boldsymbol P}^{\rm (res)}(\beta)$ are defined as the `residual counterparts' of ${\boldsymbol P}_j(\beta)$ and ${\boldsymbol P}(\beta)$, respectively (i.e., in their definitions the probabilities $P_{ij}(\beta)$ are to be replaced by $P^{\rm (res)}_{ij}(\beta)$).
\begin{lemma}\label{L2}
     For any $\beta>0$,
    \[{\boldsymbol P}^{\rm (res)}(\beta) =  {\rm diag}\{{\boldsymbol{\mathscr G}^{\rm (res)}}(\beta)\}\,Q\,{\boldsymbol P}(\beta)+\big(I- {\rm diag}\{{\boldsymbol{\mathscr G}^{\rm (res)}}(\beta)\}\big),\]
    with ${\boldsymbol P}(\beta)$ being given by Lemma \ref{L1}. 
    In addition,
    \[{\boldsymbol P}^{\rm (res)}_{0}(\beta) = {\bs 1}- {\boldsymbol P}^{\rm (res)}(\beta) {\bs 1}. \]
\end{lemma}

\subsection{Computation of the cross-moments}
The next step is to evaluate the conditional expectation $\mathbb{E}(M_i(T_{\beta})\,|\,{\boldsymbol M}(0) = {\boldsymbol m})$, which will later allow us to set up an expression for the `mixed expectation' $\E[M_j(0)M_i(T_\beta)]$.
The expression for this conditional expectation is given by, as before exploiting that the system is in stationarity at time $0$, for any $i\in\{1,\ldots,n\}$,
\[\mathbb{E}(M_i(T_{\beta})\,|\,{\boldsymbol M}(0) = {\boldsymbol m}) = \sum_{j=1}^n m_j\, P_{ji}^{\rm (res)}(\beta) + \sum_{j=1}^n \Gamma_{ji}(\beta)\]
where $\Gamma_{ji}(\beta)$ is the expected value of the number of customers that arrive at station $j$ between time~$0$ and time $T_\beta$ and that reside at station $i$ at time $T_\beta$; here it has been used that customers move through the network without interfering with each other. 

Now note that, in case the exponential clock $T_\beta$ attains the value $t>0$, then the expected number of external arrivals at station $j$ is $\lambda_j t$, each of them arriving at a uniformly distributed time in the interval $[0,t].$ Hence, with ${\mathfrak p}_{ji}(t)$ denoting the probability that a customer that  arrives at time $0$ at station $j$ is at station $i$ at time $t$,
\begin{align*}
    \Gamma_{ji}(\beta) &= \int_0^\infty \beta e^{-\beta t}  \lambda_j t\int_0^t \frac{1}{t}\,{\mathfrak p}_{ji} (t-s)\,{\rm d}s\,{\rm d}t=\lambda_j\int_0^\infty \beta e^{-\beta t} \int_0^t {\mathfrak p}_{ji} (s)\,{\rm d}s\,{\rm d}t \\
    &= \lambda_j\int_0^\infty \left(\int_s^{\infty}\beta e^{-\beta t} \,{\rm d}t\right)
     {\mathfrak p}_{ji} (s)\,{\rm d}s=\lambda_j\int_0^\infty  e^{-\beta s} \,
     {\mathfrak p}_{ji} (s)\,{\rm d}s\\
     &= \lambda_j \,\frac{P_{ji}(\beta)}{\beta}.
    \end{align*}
This can be rewritten in more compact, self-evident notation: ${\bs \Gamma}(\beta)=\beta^{-1}{\rm diag}\{{\bs\lambda}\}\,{\bs P}(\beta).$ As an aside we mention that \michel{$\sum_{j=1}^n\Gamma_{ji}(\beta)\to \varrho_i$ as $\beta\downarrow 0$, as anticipated (or ${\bs 1}^\top\Gamma(\beta)\to{\bs\varrho}^\top$ in compact notation).}

% Solving for $\Gamma_{ji}(\beta)$,

% Method 1 - 
% \begin{align*}
%     \Gamma_{ji}(\beta) &= \int_0^\infty \int_0^t \lambda_j \,\beta \,e^{-\beta t} \,g_j(x) \,P_{ji} (t-x)\,{\rm d}x\,{\rm d}t\\
%     &= \int_0^\infty \int_0^t \lambda_j \,\beta\,e^{-\beta x} \,e^{-\beta (t-x)} \,g_j(x) \,P_{ji} (t-x)\,{\rm d}x\,{\rm d}t\\
%     &= \int_0^\infty \lambda_j \,e^{-\beta x} \, \,g_j(x)\int_x^\infty \beta\,e^{-\beta (t-x)} \,P_{ji} (t-x)\,{\rm d}t\,{\rm d}x\\
%     &= \int_0^\infty \lambda_j \,e^{-\beta x} \, \,g_j(x)\int_0^\infty \beta\,e^{-\beta (t)} \,P_{ji} (t)\,{\rm d}t\,{\rm d}x\\
%     &= \int_0^\infty \lambda_j \,e^{-\beta x} \,g_j(x)\, P_{ji}(\beta)\,{\rm d}x\\
%     &= \lambda_j \,\mathscr{G}_j(\beta)\, P_{ji}(\beta)\,
% \end{align*}

% Method 2 - 

% \begin{align*}
%     \Gamma_{ji}(\beta) &= \int_0^\infty \int_0^t \lambda_j \,t \,\beta \,e^{-\beta t}\,\frac{1}{t} \,P_{ji} (t-x)\,{\rm d}x\,{\rm d}t\\
%     &= \int_0^\infty \int_0^t \lambda_j \,\beta\,e^{-\beta x} \,e^{-\beta (t-x)} \,P_{ji} (t-x)\,{\rm d}x\,{\rm d}t\\
%     &= \int_0^\infty \lambda_j \,e^{-\beta x} \,\int_x^\infty \beta\,e^{-\beta (t-x)} \,P_{ji} (t-x)\,{\rm d}t\,{\rm d}x\\
%     &= \int_0^\infty \lambda_j \,e^{-\beta x} \,\int_0^\infty \beta\,e^{-\beta (t)} \,P_{ji} (t)\,{\rm d}t\,{\rm d}x\\
%     &= \int_0^\infty \lambda_j \,e^{-\beta x} \, P_{ji}(\beta)\,{\rm d}x\\
%     &= \frac{\lambda_j \, P_{ji}(\beta)}{\beta}
% \end{align*}

Our next goal is to identify, for any $i,j\in\{1,\ldots,n\}$, the `mixed expectation' $\E[M_j(0)M_i(T_\beta)]$. By conditioning on ${\boldsymbol M}(0)={\boldsymbol m}$,
\begin{align*}\E[M_j(0)M_i(T_\beta)] &= \sum_{{\boldsymbol m}\in{\mathbb N}_0^n} \E[M_j(0)M_i(T_\beta)\,|\,{\boldsymbol M}(0)={\boldsymbol m}]\,{\mathbb P}({\boldsymbol M}(0)={\boldsymbol m})\\
&= \sum_{{\boldsymbol m}\in{\mathbb N}_0^n} m_j\E[M_i(T_\beta)\,|\,{\boldsymbol M}(0)={\boldsymbol m}]\,
\prod_{\ell = 1}^n \frac{e^{-\varrho_\ell} \,(\varrho_\ell )^{m_\ell}}{m_\ell !}\\
&=\sum_{{\boldsymbol m}\in{\mathbb N}_0^n} m_j\left(\sum_{k=1}^n m_k\, P_{ki}^{\rm (res)}(\beta) + \sum_{k=1}^n \Gamma_{ki}(\beta)
\right)\,
\prod_{\ell = 1}^n \frac{e^{-\varrho_\ell} \,(\varrho_\ell )^{m_\ell}}{m_\ell !}.\end{align*}
We thus obtain, by distinguishing the cases $j=k$ and $j\neq k$, and recalling that $\varrho^2+\varrho$ is the second moment of a Poisson random variable with mean $\varrho$, the following result:
\begin{align*}
 \E[M_j(0)M_i(T_\beta)] &= (\varrho_j^2 +\varrho_j) P_{ji}^{\rm (res)}(\beta)+\sum_{k=1,k\not=j}^n \varrho_j\varrho_k \,P_{ki}^{\rm (res)}(\beta) + \varrho_j \sum_{k=1}^n \Gamma_{ki}(\beta)\\
  &= \varrho_j \, P_{ji}^{\rm (res)}(\beta)+\varrho_j\sum_{k=1}^n \varrho_k \,P_{ki}^{\rm (res)}(\beta) + \varrho_j \sum_{k=1}^n \Gamma_{ki}(\beta).
 \end{align*}
By plugging in the relation ${\bs \Gamma}(\beta)=\beta^{-1}{\rm diag}\{{\bs\lambda}\}\,{\bs P}(\beta)$ that we derived above, in vector/matrix notation we find the following representation.

\begin{proposition} \label{P1}
For any $\beta>0$,
\begin{equation}\label{eq:cross_M}    
\E[{\bs M}(0){\bs M}(T_\beta)^\top] =\: \big( {\bs \varrho} {\bs \varrho}^\top +{\rm diag}\{{\bs \varrho}\}\big) \,{\bs P}^{\rm (res)}(\beta)+ {\beta}^{-1}\,{\bs \varrho}\,{\bs\lambda}^\top{\bs P}(\beta)
,\end{equation}
 with ${\boldsymbol P}(\beta)$ being given by Lemma \ref{L1} and ${\boldsymbol P}^{\rm (res)}(\beta)$ being given by Lemma \ref{L2}.
\end{proposition}
As a sanity check, we now verify that $\E[{\bs M}(0){\bs M}(T_\beta)^\top]\to {\bs\varrho}{\bs\varrho}^\top$ as $\beta\downarrow 0$; this property should hold, since $T_\beta\to\infty$ as $\beta\downarrow 0$, so that we should have that ${\bs M}(0)$ and ${\bs M}(T_\beta)$ become independent in the limit. \michel{To see that we indeed have that $\E[{\bs M}(0){\bs M}(T_\beta)^\top]\to {\bs\varrho}{\bs\varrho}^\top$, use that ${\bs P}^{\rm (res)}(\beta)\to 0$ and ${\bs 1}^\top\Gamma(\beta)\to{\bs\varrho}^\top$:}
\begin{align*}\lim_{\beta\downarrow 0}\E[{\bs M}(0){\bs M}(T_\beta)^\top]&=\lim_{\beta\downarrow 0}{\beta}^{-1}\,{\bs \varrho}\,{\bs\lambda}^\top{\bs P}(\beta)=\lim_{\beta\downarrow 0}{\bs\varrho}\big({\bs 1}^\top \,\beta^{-1}{\rm diag}\{{\bs\lambda}\}\,{\bs P}(\beta)\big)\\&= \lim_{\beta\downarrow 0}{\bs\varrho}{\bs 1}^\top {\bs \Gamma}(\beta)={\bs \varrho} {\bs \varrho}^\top.\end{align*}

We also remark that the result found in Proposition \ref{P1} is consistent with (known) expressions for the case of a single infinite-server queue; let the arrival rate be $\lambda$ and $G$ be the service time with Laplace-Stieltjes transform ${\mathscr G}(\cdot)$. Then $\varrho = \lambda \,\E[G]$, so that 
\[\E[{M}_1(0){M}_1(T_\beta)]= (\varrho^2+\varrho) \,{\mathbb P} (G^{\rm (res)} \geqslant T_\beta) + \varrho \int_0^\infty \beta e^{-\beta t}\,\lambda t\int_0^t \frac{1}{t}{\mathbb P}(G\geqslant t-s)\,{\rm d}s\,{\rm d}t,\]
the first term taking care of the contribution to $M_1(T_\beta)$ due to customers already present at time $0$, and the second term of the contribution  due to customers arriving between times $0$ and $T_\beta$.
After some elementary calculations, it turns out that this can be rewritten as 
\begin{equation}\label{eq:MGinff}(\varrho^2+\varrho)\,{\mathbb P} (G^{\rm (res)} \geqslant T_\beta)+\beta^{-1}\,\varrho \lambda \,{\mathbb P}(G\geqslant T_\beta),\end{equation}
as desired. It can be further simplified by expressing the probabilities ${\mathbb P} (G^{\rm (res)} \geqslant T_\beta)$ and ${\mathbb P} (G \geqslant T_\beta)$ in terms of ${\mathscr G}(\beta).$

\subsection{Erlang-distributed inter-observation times}
So far we have considered cross-moments $\E[{\bs M}(0){\bs M}(T_\beta)^\top]$ pertaining to exponentially distributed inter-observation times. In other words, we have considered the network population, in stationarity, at the $k$-th and $(k+1)$-st Poisson observation times. 
We conclude this section by pointing out how one can extend that procedure to Erlang-$2$ inter-observation times, corresponding to the cross-moments at the $k$-th and $(k+2)$-nd Poisson observation times. 

Let $E_{\beta,2}$ be an Erlang-2 random variable with mean $2\,\beta^{-1}$, i.e., the sum of two independent exponentially distributed phases, each having mean $\beta^{-1}$. The key insight is the following relation, by which one can express $\E[{\bs M}(0){\bs M}(E_{\beta,2})^\top]$ in terms of $\E[{\bs M}(0){\bs M}(T_\beta)^\top]$:
\begin{align*}
    \frac{{\rm d}}{{\rm d}\beta} \E[{\bs M}(0){\bs M}(T_\beta)^\top] &=
\frac{{\rm d}}{{\rm d}\beta} \int_0^\infty \beta e^{-\beta t}\, \E[{\bs M}(0){\bs M}(t)^\top]{\rm d}t\\
    &=\int_0^\infty  e^{-\beta t}\, \E[{\bs M}(0){\bs M}(t)^\top]{\rm d}t-\int_0^\infty \beta t e^{-\beta t}\, \E[{\bs M}(0){\bs M}(t)^\top]{\rm d}t\\
    &=\beta^{-1} \E[{\bs M}(0){\bs M}(T_\beta)^\top] - \beta^{-1} \E[{\bs M}(0){\bs M}(E_{\beta,2})^\top],
\end{align*}
where we have used that the probability density function of $E_{\beta,2}$ reads $\beta^2t\,e^{-\beta t}$, for $t\geqslant 0$. 
By rearranging this identity, we readily obtain
\begin{equation}
    \label{relErl}
\E[{\bs M}(0){\bs M}(E_{\beta,2})^\top]= \E[{\bs M}(0){\bs M}(T_\beta)^\top]-\beta\,  \frac{{\rm d}}{{\rm d}\beta} \E[{\bs M}(0){\bs M}(T_\beta)^\top].\end{equation}
In our analysis, a crucial role is played by the matrices $\dot{\bs P}(\beta)$ and $\dot{\bs P}^{\rm (res)}(\beta)$, denoting the entry-wise derivatives of ${\bs P}(\beta)$ and ${\bs P}^{\rm (res)}(\beta)$, respectively. These can be identified as follows. From \eqref{E1} we find, using the standard differentiation rules for matrix products, that 
\[(I-{\rm diag}\{{\boldsymbol{\mathscr G}}(\beta)\}Q)\,\dot{\bs P}(\beta) - {\rm diag}\{\dot{\boldsymbol{\mathscr G}}(\beta)\}Q\,{\bs P}(\beta) = -{\rm diag}\{\dot{\boldsymbol{\mathscr G}}(\beta)\},\]
so that
\begin{equation}\label{pdot}\dot{\bs P}(\beta) =(I-{\rm diag}\{{\boldsymbol{\mathscr G}}(\beta)\}Q)^{-1} \, {\rm diag}\{\dot{\boldsymbol{\mathscr G}}(\beta)\}\, (Q\,{\bs P}(\beta)-I) .\end{equation}
From this relation and Lemma \ref{L2} we can identify $\dot{\bs P}^{\rm (res)}(\beta)$:
\begin{equation}\label{pdotres}\dot{\bs P}^{\rm (res)}(\beta) ={\rm diag}\{\dot{\boldsymbol{\mathscr G}}^{\rm (res)}(\beta)\}\,Q\, {\bs P}(\beta)+{\rm diag}\{{\boldsymbol{\mathscr G}}^{\rm (res)}(\beta)\}\,Q\, \dot {\bs P}(\beta)-
{\rm diag}\{\dot{\boldsymbol{\mathscr G}}^{\rm (res)}(\beta)\}.
\end{equation}
Applying \eqref{relErl} in combination with \michel{Proposition \ref{P1}}, we thus end up with the following result.
\begin{proposition}\label{P2}
For any $\beta>0$,
\begin{align*}\E[{\bs M}(0){\bs M}(E_{\beta,2})^\top]=&\:\big( {\bs \varrho} {\bs \varrho}^\top +{\rm diag}\{{\bs \varrho}\}\big) \,\big({\bs P}^{\rm (res)}(\beta)-\beta \dot{\bs P}^{\rm (res)}(\beta)\big)
%\:+\\&\:
+2\,{\beta}^{-1}\,{\bs \varrho}\,{\bs\lambda}^\top{\bs P}(\beta)-{\bs \varrho}\,{\bs\lambda}^\top\dot{\bs P}(\beta),
\end{align*}
with $\dot{\bs P}(\beta)$ being given by \eqref{pdot} and $\dot{\bs P}^{\rm (res)}(\beta)$ by \eqref{pdotres}.
\end{proposition}
Along the same lines, the counterpart of Proposition \ref{P2} for higher-order Erlang inter-observation times can be found. One more differentiation and similar manipulations, for instance, lead us to an expression for $\E[{\bs M}(0){\bs M}(E_{\beta,3})^\top].$

%%%
\section{Method of moments estimator and its consistency}\label{sec:method}

As pointed out in the preceding sections, the parameters of our interest are the routing matrix $Q$, the arrival rates ${\bs\lambda}$, and the parameters associated to the service-time distributions. In this section we assume that the service-time distributions  are given, leaving us with the task of estimating the $n$ arrival rates in the vector ${\bs\lambda}$, the $n^2$ entries of the routing matrix ${Q}$, yielding a total of $n^2+n$ parameters. As mentioned in the introduction, Section \ref{sec:G_unknown} will address how inference of the service-time distributions can be incorporated, while Section \ref{sec:obsprob} will discuss how to handle noisy observations, which requires estimating the 
$n$ observation probabilities as well.

The main objective of this section is to propose and analyze an estimator for the parameters $Q$ and~$\bs\lambda$. The moment equations we will work with are those related to $\E[{\bs M}(0)]$ and $\E[{\bs M}(0){\bs M}(T_\beta)^\top]$, as given in Section \ref{sec:mom}. These equations total $n^2+n$, matching the number of parameters to be estimated. Our primary result concerns the consistency of the estimator of $Q$ and $\bs\lambda$, under standard regularity assumptions.

\subsection{Estimator and assumptions} \label{ssec:est}In this subsection we define our estimator, and state the assumptions under which we can prove consistency. 
Throughout we let our parameter vector be represented through
\begin{align*}
   \michel{ {\bs \theta}\equiv (Q,{\bs\lambda})\in\Theta,\qquad \Theta:=\left\lbrace Q\in \Delta_{n}, {\bs \lambda}\in[0,\Lambda]^n\   \right\rbrace };
\end{align*}
where $\Delta_{n}$ consists of all $n\times n$ matrices with non-negative entries in which each row sum is less than or equal to one., and where $\Lambda\in [0, \infty)$ is to be interpreted as an upper bound on each of the external input rates.  \michel{We denote the true parameter by ${\bs \theta}_0\equiv ( Q_0,{\bs \lambda}_{0})$, being an element of the parameter space $\Theta$.}

Let $\tau_1,\ldots,\tau_m$ denote the sequence of inter-sampling times, with $m\in{\mathbb N}$ denoting the number of observations. As we are working with Poisson observation times, these inter-sampling times are independent, exponentially distributed samples, say with mean $\beta^{-1}$; the sampling rate $\beta$ is omitted from our notation for the sake of brevity. The sampling times are then given by the corresponding partial sums: $S_k=\sum_{j=1}^k\tau_i$, with $k\in\{1,\ldots,m\}$. In the sequel, we denote the queue size at station $i$ and observation time $k$ by $M_{i,k}:=M_i(S_k)$, and the vector recording the queue sizes at all $n$ stations by ${\bs M}_k$. %Furthermore, denote the these quantities' {\it observed} counterparts by $N_{i,k}:=N_i(S_k)$ and ${\bs N}_k$. 

The next step is to use the expressions found in Section \ref{sec:mom} to set up the moment equations that we use in our estimator. 
To this end, for any ${\bs\theta}\in\Theta$, we define two types of random objects:
\begin{align*}
    \psi_{i,k}^{[0]}(\bs\theta)&:=M_{i,k}-\,{\mathbb E}_{\bs\theta}[M_i(0)]=M_{i,k}-\,\lambda_i^{\mathrm{eff}}\,\mathbb{E}[G_i]\ ,
    \end{align*}
for $i\in\{1,\ldots,n\}$ and $k\in\{1,\ldots,m\}$, and
    \begin{align*}
    \psi_{i,j,k}^{[1]}(\bs\theta)&:=M_{j,k}M_{i,k+1}-\,\mathbb{E}_{\bs\theta}[M_j(0)M_i(T_\beta)]\ ,
\end{align*}
for $i,j\in\{1,\ldots,n\}$ and $k\in\{1,\ldots,m-1\}$. Observe  that all these quantities have mean zero for $\bs\theta=\bs\theta_0$. \liron{The method-of-moments estimator is based on solving the following system of moment equations:
\begin{align}\label{ME1}
    \frac{1}{m}\sum_{k=1}^m\psi^{[0]}_{i,k}(\tilde{\bs \theta}_m) = 0,&\quad i\in\{1,\ldots,n\}\ , \\ \label{ME2}
      \frac{1}{m-1}\sum_{k=1}^{m-1}\psi^{[1]}_{i,j,k}(\tilde{\bs \theta}_m) = 0,& \quad i,j\in\{1,\ldots,n \} \ .
\end{align}
It is noted, however, that the solution $\tilde{\bs\theta}_m$ might yield estimators that are outside the parameter space $\Theta$. To remedy this, we follow the standard procedure of working with its projected counterpart $\hat{\bs\theta}_m$, i.e., we set, with $x_+ := \max\{x,0\}$ denoting the positive part of $x$,
\begin{align}
    \hat{\lambda}_{i,m}
    &:= \min\!\left\{(\tilde{\lambda}_{i,m})_+,\Lambda\right\},
    \quad i=1,\ldots,n, \label{eq:lambda_truncated}\\[2pt]
    \hat{q}_{i,j,m}
    &:= \frac{(\tilde{q}_{i,j,m})_+}
    {\max\!\bigg\{1,\displaystyle\sum_{k=1}^n(\tilde{q}_{i,k,m})_+\bigg\}},
    \quad i,j=1,\ldots,n. \label{eq:q:truncated}
\end{align}

Observe that by construction  $\hat{\bs\theta}_m\in\Theta$. Throughout the following assumptions are in place. }

\begin{assumption}\label{assum:Q_flow}
\liron{The true parameter $\bs\theta_0= (Q,{\bs \lambda})$ satisfies \eqref{eq:C1} and \eqref{eq:C2}.
% \begin{enumerate}
% %
%     \item[(i)] For any station $j$,
%     \begin{equation}
% \notag     \lambda_j>0, \:\: \mbox{or there is a path $(i_0, \ldots,i_k)$ such that $\lambda_{i_0}q_{i_0,i_1}\cdots q_{i_k j}>0$.}
%   \end{equation} 
%     %
%     \item[(ii)] There exists a positive constant $\epsilon$ such that for any station $j=1,\ldots,n$,
%     \begin{equation}\notag
%    q_{j0}>\epsilon,\:\: \mbox{or there is a path $(i_0, \ldots,i_k)$ such that $q_{ji_0}\cdots q_{i_{k-1} i_k}q_{i_k0}>\epsilon$.}   
%   \end{equation}
% \end{enumerate}
    }
\end{assumption}

% In our setup we impose the following two conditions on the model parameters. (i)~In the first place we want each queue $j$ to be `non-trivial' in that customers flow in. This is guaranteed if \begin{equation}\label{eq:C1}
%      \lambda_j>0, \:\: \mbox{or there is a path $(i_0, \ldots,i_k)$ such that $\lambda_{i_0}q_{i_0,i_1}\cdots q_{i_k j}>0$.}
%  \end{equation} (ii)~In the second place we want each queue $j$ to be `non-explosive', in the sense that customers eventually flow out. This is guaranteed if 
%  \begin{equation}\label{eq:C2}
%   q_{j0}>0,\:\: \mbox{or there is a path $(i_0, \ldots,i_k)$ such that $q_{ji_0}\cdots q_{i_{k-1} i_k}q_{i_k0}>0$.}   
%  \end{equation}

%%%
\begin{assumption}\label{assum:G_moments}
For any $j=1,\ldots,n$, $\mathbb{E}[G_j]$ is positive and finite.
\end{assumption}

We proceed by providing an intuitive explanation as to why these assumptions need to be imposed. 
Assumption~\ref{assum:Q_flow} specifies that all nodes have a non-zero rate of incoming work and that all jobs leave the system eventually. 
%These conditions are similar to the ones imposed in \eqref{eq:C1}--\eqref{eq:C2}, but now it should be in place for all ${\bs\theta}\in\Theta$, where it is noted that in part (ii) we need the uniform lower bound $\epsilon$ to make the consistency proof work. 
\liron{Assumption~\ref{assum:G_moments} ensures, together with Assumption~\ref{assum:Q_flow}, system stability and almost sure convergence of the empirical moments \eqref{ME1}--\eqref{ME2}, as appearing in the moment equations.} The assumption that all mean service-time times are strictly positive excludes the degenerate case where all jobs in a certain node have zero service time and immediately depart; this is without any loss of generality, as  such a node can be removed from the network.

\subsection{Identifiability}

In this subsection we wish to identify the unique solution to the moment equations, recalling that we focus on the case that the service-time distributions are known. A statistical model is identifiable if $\mathbb{P}_\theta=\mathbb{P}_{\tilde{\theta}}$ if and only if $\theta=\tilde{\theta}$, where $\mathbb{P}_\theta$ is the joint distribution of the data corresponding to parameter $\theta$ (see e.g., \cite[Ch.~5]{V2000}). 
This is not sufficient for consistency of method-of-moment estimators, but rather we need to verify that there is a one-to-one mapping between the moments and the parameters.
We next verify that this is indeed the case in our setting.

For brevity, denote by the vector
$\bar{\bs \alpha}^{[0]}$ the stationary moments ${\mathbb E}[{\bs M(0)}]$ and by the matrix $\bar{\bs\alpha}^{[1]}$  the stationary cross-moments of ${\mathbb E}[{\bs M(0)}{\bs M}(T_\beta)^\top]$. Firstly, \eqref{eq:loads} and \eqref{eq:cross_M} imply that any parameter pair $(Q,\bs\lambda)$ yields a unique pair $(\bar{\bs \alpha}^{[0]},\bar{\bs \alpha}^{[1]})$. We now verify that the other direction holds as well. 

Evidently, from \eqref{eq:loads} and \eqref{eq:loads2}, we find, for $i\in\{1,\ldots,n\}$,
\begin{equation}\label{eq:estlambda}\lambda_i^{\rm eff}(\bar{\bs \alpha}^{[0]}) = \frac{\bar\alpha_i^{[0]}}{{\mathbb E}[G_i]}.\end{equation}
Now consider the moment equation, with the second equality being due to \eqref{eq:leff},
\begin{align*}\bar{\bs\alpha}^{[1]} &=\big( {\bs \varrho} {\bs \varrho}^\top +{\rm diag}\{{\bs \varrho}\}\big) \,{\bs P}^{\rm (res)}(\beta)+ {\beta}^{-1}\,{\bs \varrho}\,{\bs\lambda}^\top{\bs P}(\beta),\\
&=\big( {\bs \varrho} {\bs \varrho}^\top +{\rm diag}\{{\bs \varrho}\}\big) \,{\bs P}^{\rm (res)}(\beta)+ {\beta}^{-1}\,{\bs \varrho}\,({\bs\lambda}^{\rm eff}(\bar{\bs \alpha}^{[0]}))^\top(I-Q){\bs P}(\beta).
\end{align*}
In this equation we replace $\bs\varrho$ by the moment vector $\bar{\bs\alpha}^{[0]}$, and likewise ${\bs\lambda}^{\rm eff}$ by ${\bs\lambda}^{\rm eff}(\bar{\bs \alpha}^{[0]})$.
By Lemmas~\ref{L1} and~\ref{L2} we have expressions for the matrices ${\bs P}(\beta)$ and ${\bs P}^{\rm (res)}(\beta)$. Inserting these into the previous display, we readily arrive at 
\begin{align*}
\bar{\bs\alpha}^{[1]} =\:& \Phi\, \Big({\rm diag}\{{\boldsymbol{\mathscr G}^{\rm (res)}}(\beta)\}Q\,\big((I-{\rm diag}\{{\boldsymbol{\mathscr G}(\beta)\}Q})^{-1}(I-{\rm diag}\{{\boldsymbol{\mathscr G}}(\beta)\})\big)+I-{\rm diag}\{{\boldsymbol{\mathscr G}^{\michel{\rm (res)}}}(\beta)\}\Big)\,+\\
& {\beta}^{-1}\,\bar{\bs\alpha}^{[0]}\,(\bs\lambda^{\rm eff}(\bar{\bs \alpha}^{[0]}))^\top(I-Q)\,\big((I-{\rm diag}\{{\boldsymbol{\mathscr G}(\beta)\}Q})^{-1}(I-{\rm diag}\{{\boldsymbol{\mathscr G}}(\beta)\})\big),
\end{align*}
where
$\Phi := \bar{\bs\alpha}^{[0]} (\bar{\bs\alpha}^{[0]})^\top+{\rm diag}\{\bar{\bs\alpha}^{[0]}\}.$
It remains to show how $Q$ can be isolated from this equation. The idea is to post-multiply the entire equation by $(I-{\rm diag}\{{\boldsymbol{\mathscr G}}(\beta)\})^{-1}(I-{\rm diag}\{{\boldsymbol{\mathscr G}(\beta)\}Q})$; for compactness we rewrite this matrix as 
$\Xi_1 +\Xi_2\,Q,$
where 
\[\Xi_1\equiv\Xi_1(\beta):=(I-{\rm diag}\{{\boldsymbol{\mathscr G}}(\beta)\})^{-1},\:\:\:\:\Xi_2\equiv \Xi_2(\beta):=-(I-{\rm diag}\{{\boldsymbol{\mathscr G}}(\beta)\})^{-1}{\rm diag}\{{\boldsymbol{\mathscr G}}(\beta)\}.\]
After elementary manipulations, we thus obtain the equation
\begin{align*}
    \bar{\bs\alpha}^{[1]}&\,\big(\Xi_1 +\Xi_2\,Q\big) \\=\:&\Phi\, \Big({\rm diag}\{{\boldsymbol{\mathscr G}^{\rm (res)}}(\beta)\}Q+\big(I-{\rm diag}\{{\boldsymbol{\mathscr G}^{\michel{\rm (res)}}}(\beta)\}\big)\big(\Xi_1 +\Xi_2\,Q\big)\Big)+{\beta}^{-1}\,\bar{\bs\alpha}^{[0]}\,({\bs\lambda}^{\rm eff}(\bar{\bs \alpha}^{[0]}))^\top(I-Q).
\end{align*}
Observe that this equation is linear in the matrix $Q$. \michel{More precisely,
with the matrices $\Omega_1\equiv \Omega_1 (\bar{\bs \alpha}^{[0]},\bar{\bs \alpha}^{[1]})$ and $\Omega_2\equiv \Omega_2 (\bar{\bs \alpha}^{[0]},\bar{\bs \alpha}^{[1]})$ given by
\begin{align*}
    \Omega_1 &:=\bar{\bs\alpha}^{[1]}\,\Xi_2 -\Phi\Big({\rm diag}\{{\boldsymbol{\mathscr G}^{\rm (res)}}(\beta)\}+\big(I-{\rm diag}\{{\boldsymbol{\mathscr G}^{\michel{\rm (res)}}}(\beta)\}\big)\Xi_2\Big) + {\beta}^{-1}\,\bar{\bs\alpha}^{[0]}\,({\bs\lambda}^{\rm eff}(\bar{\bs \alpha}^{[0]}))^\top,\\
    \Omega_2&:=-\bar{\bs\alpha}^{[1]}\,\Xi_1 +\Phi\,\big(I-{\rm diag}\{{\boldsymbol{\mathscr G}^{\michel{\rm (res)}}}(\beta)\}\big)\Xi_1+{\beta}^{-1}\,\bar{\bs\alpha}^{[0]}\,({\bs\lambda}^{\rm eff}(\bar{\bs \alpha}^{[0]}))^\top,
\end{align*}
respectively,
we have thus identified the mapping of moments to the pair $(Q,\bs\lambda)$. }

\begin{proposition}\label{prop:identify}
     Let ${\bs\lambda}^{\rm eff}(\bar{\bs \alpha}^{[0]})$ be given by \eqref{eq:estlambda}. Then, \michel{provided that $\Omega_1(\bar{\bs \alpha}^{[0]},\bar{\bs \alpha}^{[1]})$ is invertible,}
    \begin{align*} Q(\bar{\bs \alpha}^{[0]},\bar{\bs \alpha}^{[1]})&= \Omega_1(\bar{\bs \alpha}^{[0]},\bar{\bs \alpha}^{[1]})^{-1}\,\Omega_2(\bar{\bs \alpha}^{[0]},\bar{\bs \alpha}^{[1]}),\\
    {\boldsymbol\lambda}(\bar{\bs \alpha}^{[0]}) &= (I- Q(\bar{\bs \alpha}^{[0]},\bar{\bs \alpha}^{[1]})^\top)\, {\boldsymbol \lambda}^{\rm eff}(\bar{\bs \alpha}^{[0]}).\end{align*}
\end{proposition}
A direct corollary of the proposition is that there is a one-to-one relation between the moments $(\bar{\bs\alpha}^{[0]},\bar{\bs\alpha}^{[1]})$ and the parameters $( Q,{\bs\lambda})$. 

\begin{remark}\label{rem:R1}\em
  Proposition~\ref{prop:identify} further has practical implications because the parameter estimators can be computed by replacing $(\bar{\bs \alpha}^{[0]},\bar{\bs \alpha}^{[1]})$ with their empirical counterparts. Note, however, that for finite samples this may yield estimators outside the parameter space \liron{and projection back to the parameter space is therefore required; \michel{see \eqref{eq:lambda_truncated}--\eqref{eq:q:truncated}}. In Section~\ref{sec:exp}, we provide further details on the practical implementation of the estimation procedure.}
  \hfill $\Diamond$ 
\end{remark}
%  , for instance, negative arrival rates or routing matrices that are not stochastic. In such cases, some truncation or projection onto the admissible parameter space may be required. 

\begin{remark}\label{rem:R1a}\em
\michel{Another practical issue concerns the invertibility of the matrix $\Omega_1$, which is required for the explicit expressions in Proposition~\ref{prop:identify}. Failure of invertibility can only occur for a very specific set of values of the estimated counterparts of $(\bar{\bs \alpha}^{[0]},\bar{\bs \alpha}^{[1]})$ and a very specific sampling rate $\beta$. More precisely, singularity is characterized by the condition $\det(\Omega_1(\bar{\bs \alpha}^{[0]},\bar{\bs \alpha}^{[1]}))=0$, which defines an `exceptional subset' of the corresponding moment space, in that under natural regularity conditions it is negligible (in particular, it has Lebesgue measure zero). Consequently, non-invertibility arises only in highly degenerate situations and is not expected to affect the practical applicability of the proposed estimation procedure. }\hfill $\Diamond$ 
\end{remark}

\begin{remark}\hritika{{\em When the system is sufficiently simple, we can derive explicit expressions for $({\bs\lambda},Q)$ in terms of the first moments and cross moments. We illustrate this in the single-node setting, where the moments $\bar{\bs\alpha}^{[0]}$ and $\bar{\bs\alpha}^{[1]}$ reduce to scalars:
\[
\bar\alpha^{[0]}=\E[M_1(0)], \qquad
\bar\alpha^{[1]}=\E[M_1(0)M_1(T_\beta)].
\]
For instance in the case of exponential service times (entailing that $\mathscr G^{\rm(res)}(\beta)=\mathscr G(\beta)$), Proposition~\ref{prop:identify} yields the equations
\begin{align*}
q_{11}(\bar\alpha^{[0]},\bar\alpha^{[1]})
&=1+\frac{\beta}{\mu_1}
-\frac{\beta\,\bar\alpha^{[0]}}
{\mu_1\big(\bar\alpha^{[1]}-(\bar\alpha^{[0]})^2\big)}
=1+\frac{\beta\big(\bar\alpha^{[1]}-\bar\alpha^{[0]}-(\bar\alpha^{[0]})^2\big)}
{\mu_1\big(\bar\alpha^{[1]}-(\bar\alpha^{[0]})^2\big)},\\
\lambda_1(\bar\alpha^{[0]},\bar\alpha^{[1]})
&=\frac{\beta(\bar\alpha^{[0]})^2}
{\bar\alpha^{[1]}-(\bar\alpha^{[0]})^2}
-\beta\bar\alpha^{[0]}.
\end{align*}
It is readily verified that $q_{11}\in[0,1]$ and $\lambda_1\geqslant 0$. 
For instance, $q_{11}\leqslant 1$ can be shown as follows. Since the marginal distribution of $M_1(0)$ has a Poisson distribution,
\[
\Var[M_1(0)]
=\E[M_1(0)]
=\bar\alpha^{[0]}.
\]
Now it follows that $\bar\alpha^{[1]}-\bar\alpha^{[0]}-(\bar\alpha^{[0]})^2\leqslant 0$ by observing
\[
\bar\alpha^{[1]}-(\bar\alpha^{[0]})^2
=\E[M_1(0)M_1(T_\beta)]-\E[M_1(0)]^2
=\Cov(M_1(0),M_1(T_\beta))
\leqslant \Var[M_1(0)]
=\bar\alpha^{[0]},
\]
where the inequality follows due to the Cauchy--Schwarz inequality. 
Noting that
$
\bar\alpha^{[1]}-(\bar\alpha^{[0]})^2
=\Cov(M_1(0),M_1(T_\beta))>0,$ it is now immediate that $q_{11}\leqslant 1$.}
% Furthermore, the theoretical quantity $\alpha^{[1]}$ is a function of the resampling rate $\beta$. Substituting its explicit expression into the formulas for $Q$ and $\lambda$ shows that the dependence on $\beta$ cancels exactly, as expected. In our estimation procedure, however, $\alpha^{[1]}$ is replaced by its empirical counterpart. Consequently, the estimators of $Q$ and $\lambda$ do depend on the choice of $\beta$, where it is noted that the statistical accuracy of the estimator for $\alpha^{[1]}$ depends on the observation interval determined by $\beta$. A small value of $\beta$ results in infrequent sampling and therefore weakly correlated observations, whereas a large value of $\beta$ leads to frequent sampling and consequently strongly correlated observations. This trade-off suggests the existence of an optimal choice of $\beta$ that minimizes the estimator's variance.}
\hfill $\Diamond$}
\end{remark}

\subsection{Consistency}

% \footnote{\textcolor{blue}{Is a new subsection really needed here? This is for the same setting as 4.1, and consistency in the title of Section 4.}}

In this subsection we argue that the proposed estimator is, under the assumptions imposed,  consistent.

\liron{
%%%
\begin{theorem}\label{thm:consistency}
    Suppose that Assumptions~\ref{assum:Q_flow}--\ref{assum:G_moments} hold and that $\Omega_1(\bar{\bs \alpha}^{[0]},\bar{\bs \alpha}^{[1]})$ is invertible. Then, the estimator is consistent, i.e.,  for any sampling rate $\beta>0$, 
    \[\hat{\bs \theta}_m\asarrow  {\bs \theta}_0,\:\: \mbox{as}\:\: m\to\infty.\]
\end{theorem}
\begin{proof}
The proof of Theorem~\ref{thm:consistency} consists of three steps. The probabilistic arguments below are all with respect to the probability measure induced by $\bs\theta_0$.
    \begin{enumerate}
    \item[(i)] We employ a coupling argument to construct an (upper) bounding system.  This majorizing process is the queue length process of a specific M/G/$\infty$ queue, i.e., a queue with Poisson arrivals, i.i.d.\ service times, and infinitely many servers; the precise construction of the majorizing  M/G/$\infty$ queue follows later in Section \ref{sec:bound}. 
        \noindent The majorizing process, in the sequel denoted by $(\tilde{M}_{k})_{k\in{\mathbb N}}$, is set up such that $M_{i,k}\leqslant \tilde{M}_{k}$ for all $i\in\{1,\ldots,n\}$ and all $k\in\{1,\ldots,m\}$, with probability one. 
        %this holds for any ${\bs \theta}\in\Theta$, where we note that the majorizing process depends on ${\bs\theta}$ as well.
        We further verify that under Assumptions~\ref{assum:Q_flow}--\ref{assum:G_moments}, this majorizing process has ergodic limits for the first moment and the lag-one cross-moment. Indeed, we argue that
        %for any $\bs\theta\in\Theta$ and for 
        there exist finite numbers $\bar M$ and $\bar M_2$ such that
    \begin{equation}\label{eq:ergodic_bound}
    \frac{1}{m}\sum_{k=1}^m\tilde{M}_{k}\asarrow \bar{M}\ , \:\:\: \frac{1}{m-1}\sum_{k=1}^{m-1}\tilde{M}_{k}\tilde{M}_{k+1}\asarrow \bar{M}_2 \ ,
    \end{equation}
     as $m\to\infty$.
    \item[(ii)]  By \eqref{eq:ergodic_bound}, in combination with the fact that $(\tilde{M}_{k})_{k\in{\mathbb N}}$ majorizes $({M}_{i,k})_{k\in{\mathbb N}}$ for $i\in\{1,\ldots,n\}$, we conclude that
    \begin{align*}
         \lim_{m\to\infty}\frac{1}{m}\sum_{k=1}^m M_{i,k}<\infty \ ,
    \end{align*}
    almost surely. As we have assumed that the system is in stationarity, we have that $\E[M_{i,k}]=\lambda_i^{\mathrm{eff}}\mathbb{E}[G_i]$ for all $k\in\{1,\ldots,m\}$ and $i\in\{1,\ldots,n\}$. As a consequence, C\'esaro’s lemma implies that, for any $i\in\{1,\ldots,n\}$,
    \begin{align*}
        \frac{1}{m}\sum_{k=1}^m M_{i,k}\asarrow  \lambda_{i,0}^{\mathrm{eff}}\,\mathbb{E}[G_i] \ ,
    \end{align*}
    with ${\bs \lambda}_{0}^{\mathrm{eff}}\equiv(\lambda_{1,0}^{\mathrm{eff}},\ldots,\lambda_{n,0}^{\mathrm{eff}})^\top$ defined in the evident manner. 
    Following an identical argumentation, we can verify the limit, for any $i,j\in\{1,\ldots,n\}$,
    \begin{align*}
        \frac{1}{m-1}\sum_{k=1}^{m-1}M_{j,k}M_{i,k+1}&\asarrow   \E_{{\bs\theta}_0}\big[M_j(0)M_i(T_\beta)\big]<\infty\  ,
    \end{align*}
    again in evident notation.
    \item[(iii)]  Assuming \eqref{eq:ergodic_bound} we have that the empirical moments converge to $(\bar{\bs \alpha}^{[0]},\bar{\bs \alpha}^{[1]})$ almost surely as $m\to\infty$. If $\Omega_1(\bar{\bs \alpha}^{[0]},\bar{\bs \alpha}^{[1]})$ is invertible then Proposition~\ref{prop:identify} implies that the mapping $(\bar{\bs \alpha}^{[0]},\bar{\bs \alpha}^{[1]})\mapsto (Q,\bs\lambda)$ is continuous at the true parameter values. The truncation map \eqref{eq:lambda_truncated}-\eqref{eq:q:truncated} is continuous on all of $\mathbb{R}^{n^2+n}$ and restricts to the identity on $\Theta$, since $\bs\theta_0\in\Theta$ it is fixed by the map, whether or not $\bs\theta_0$ lies on the interior or boundary of $\Theta$. Then, we can invoke the continuous mapping theorem (see, e.g., \cite[Thm.~2.3]{V2000}) to conclude the proof.     Hence, it is left to prove that \eqref{eq:ergodic_bound} holds, which will be done in Subsection~\ref{sec:bound}.
\end{enumerate}
This completes the proof.
% Recall that the routing probabilities appearing in the parameter vector $\bs\theta$ should be component-wise nonnegative and such that $\sum_{j=1}^n q_{ij}\leqslant 1$ for any $i\in\{1,\ldots,n\}$, i.e., they belong to specific simplices. In addition, we assumed $\lambda_i\in[0,\Lambda]$ for $\Lambda<\infty$,  for all $i\in\{1,\ldots,n\}$. Hence, the parameter space $\Theta$ is compact. 
% Observe from \eqref{eq:defpsi} that $\bs\Psi_k(\bs\theta):\Theta\to\mathbb{R}^{n^2+n}$ is a vector-valued function. With the system being in stationarity at time zero, Proposition~\ref{prop:pointwise_conv} entails that
% \begin{equation}\label{eq:pointwise_conv}
%     \frac{1}{m}\sum_{k=1}^m\bs\Psi_k(\bs\theta)\asarrow \mathbb{E}[\bs\Psi_1(\bs\theta)]\ , \:\forall\theta\in\Theta\ .
% \end{equation}
% The Lipschitz continuity of $\bs\Psi_k(\cdot)$, stated in Lemma~\ref{lem:lip},  combined  with the pointwise convergence in \eqref{eq:pointwise_conv}, thus yields that $\{m^{-1}\sum_{k=1}^m\bs\Psi_k(\bs\theta): \ \bs\theta\in\Theta\ , \michel{m\in{\mathbb N}}\}$ is a Glivenko-Cantelli class (see, e.g., \cite[Thm.~3]{A1992}). Hence, as $m\to\infty$,
% \begin{align*}
%   \sup_{\bs\theta\in\Theta} \left|\frac{1}{m}\sum_{k=1}^m\bs\Psi_k(\bs\theta)-\mathbb{E}[\bs\Psi_1(\bs\theta)]\right|\asarrow 0\ .
% \end{align*}
% As $\mathbb{E}[\bs\Psi_1(\bs\theta)]=\bs{0}$ if and only if $\bs\theta=\bs\theta_0$ by Proposition~\ref{prop:pointwise_conv}, we conclude that the estimator is consistent (see, e.g., \cite[Thm.~5.9]{V2000}). 
\end{proof}
}
%%%
\subsection{Upper bound system}\label{sec:bound}
The objective of this subsection is to verify the almost sure convergences in~\eqref{eq:ergodic_bound}. In this proof we work with a coupling which features a majorizing process that is of the M/G/$\infty$ type. 

Let $M^+(t)={\bs M}(t)^\top{\bs 1}$ denote the total network population at time $t\geqslant 0$. The key observation is that $M^+(\cdot)$ is \liron{a stationary and ergodic} M/G/$\infty$ process. In the first place, its arrival process is Poissonian with rate $\lambda:=\sum_{i=1}^n\lambda_i$.
Secondly, the clients' service times are independent and identically distributed random variables, distributed as a generic random variable $H$, independent of the arrival process, that can be described as  follows. Let $H_i$ be the generic random variable that represents how long a job stays in the network, given it enters at node $i$; then \liron{$H\equiv H({\bs\theta}_0)$ equals $H_i$} with probability $\lambda_i/\lambda.$ Here it has been used that the  times that the customers spend in the network are independent, as a direct consequence of the underlying infinite-server mechanism.
We observe that, evidently, $M_i(t) \leqslant M^+(t)$ almost surely, for any $t$ and $i=1,\ldots,n.$

Let $K_i$ denote the number of nodes visited by an arbitrary customer entering the network at node $i=1,\ldots,n.$ Then, with $g:=\max\{{\mathbb E}[G_1],\ldots,{\mathbb E}[G_n]\}$, we have that ${\mathbb E}[H_i]\leqslant g\,{\mathbb E}[K_i]$; here it has been used that, by Assumption \ref{assum:G_moments}, $g$ is finite. By Assumption~\ref{assum:Q_flow} we conclude that ${\mathbb E}[K_i]$ is  finite; to this end observe that $K_i$ can be interpreted as the hitting time of an absorbing state in a discrete-time Markov chain with finite state space. \liron{This in particular means that ${\mathbb E}[H]<\infty$. }
\liron{The empirical mean of an ergodic stationary process sampled according to an independent Poisson process converges to the same stationary mean   \cite{K1999}, hence}
%By the ergodic theorem for stationary processes, recalling that $(M^+(S_k))_{k\in{\mathbb N}}$ is stationary, it follows that, as $m\to\infty$,
\begin{align*}
     \frac{1}{m}\sum_{k=1}^m M^+(S_k)\to_{\rm as}{\mathbb E}[M^+(S_1)] = \lambda\,\E[H],
 \end{align*}
 using that $M^+(S_1)\sim {\rm Poiss}(\lambda\,\E[H])$, with $\lambda\,\E[H]<\infty$.
% \begin{align*}
% \lim_{t\to\infty}    \frac{1}{t}\int_0^t M^+(u)\diff u=\E[M^+]\ ,
% \end{align*}
% where $M^+\sim {\rm Poiss}(\lambda\,\E[H])$, with $\lambda\,\E[H]<\infty$. By the {\small PASTA} property, the same is true for the empirical process at Poisson sampling times: 
% \begin{align*}
%     \frac{1}{m}\sum_{k=1}^m M^+(S_k)\to_{\rm as}\lambda\,\E[H].
% \end{align*}
We thus conclude that the first claim in \eqref{eq:ergodic_bound} holds. 
% Finally, note that by Assumption~\ref{assum:G_moments} the bounds provided here are uniform on $\theta$, i.e., there exist constants $\bar{s}_1,\bar{s}_2<\infty$ such that
% \begin{align*}
%     \sup_{\theta\in\Theta}\mathbb{E}_\theta[S] <\bar{s}_1\ ,  \sup_{\theta\in\Theta}\mathbb{E}_\theta[S^2] <\bar{s}_2\ . \end{align*}
Regarding the second claim in \eqref{eq:ergodic_bound}, we observe that $(M^+(S_k)M^+(S_{k+1}))_{k\in{\mathbb N}}$ is stationary. \liron{By the same ergodic theorem}, we now obtain, as $m\to\infty$, 
\[
 \frac{1}{m-1}\sum_{k=1}^{m-1} M^+(S_k)M^+(S_{k+1})\to {\mathbb E}[M^+(S_1)M^+(S_{2})],\]
which, in view of \eqref{eq:MGinff}, leads to, as $m\to\infty$,
\[
 \frac{1}{m-1}\sum_{k=1}^{m-1} M^+(S_k)M^+(S_{k+1})\to_{\rm as}(\varrho_+^2+\varrho_+)\,{\mathbb P} (H^{\rm (res)} \geqslant T_\beta)+\beta^{-1}\,\varrho_+ \lambda \,{\mathbb P}(H\geqslant T_\beta)<\infty,\]
 with $\varrho_+:=\lambda\,{\mathbb E}[H].$

\section{Unknown service-time distribution}\label{sec:G_unknown}

% \textcolor{blue}{Just pasted some old text here - need to edit. Is this also the place to discuss the derivative based estimation equations?}

\michel{So far, we have assumed that the service-time distributions are known. In this section, we consider the more general setting in which they are unknown. We derive the corresponding moment expressions, which lead naturally to an extension of the estimation framework in which the service-time parameters are estimated jointly with the network parameters. This extension substantially broadens the applicability of the proposed methodology, and we investigate its practical performance in the numerical experiments below. A theoretical treatment of the extended estimator would require additional identifiability conditions and a separate consistency analysis, which are beyond the scope of the present paper.
}

\subsection{Service-time distributions in parametric family} We start by discussing the setting in which the service-time distributions are in a given parametric family. Let for ease  each service-time distribution be characterized by $d$ parameters.

We let our parameter vector be represented through
\begin{align*}
    \michel{{\bs \theta}\equiv(Q,\bs\lambda,{\bs\eta)\in\Theta,\qquad \Theta:=} \left\lbrace Q\in\Delta_n, {\bs \lambda}\in[0,\Lambda]^n\  ,  {\bs \eta}\in\mathscr{E}^n  \right\rbrace,}
\end{align*}
where $\mathscr{E}\subseteq\mathbb{R}^d$
is some parametric space; it means that have that $G_j\equiv G_j(\bs\eta_j)$,
with $\bs\eta_j\in{\mathscr E}$, for any $j\in\{1,\ldots,n\}$. As before, we refer to the set of such parameters as $\Theta$, and we denote the true parameter by ${\bs \theta}_0\equiv(Q_0, {\bs \lambda}_{0},{\bs \eta}_0)\in\Theta$. 

As discussed in Remark \ref{rem:R1}, we now impose the condition $q_{ii}=0$ for all $i\in\{1,\ldots,n\}$, to avoid any identifiability issues.
Observe that the dimension of the parameter vector has now become $n^2+dn$: $n^2-n$ entries in $Q$, $n$ entries in $\bs\lambda$ and $dn$ entries in $\bs\eta$. This means that unless $d=1$ we need to add cross-moments pertaining to Erlang inter-observation times, where as long as $d\leqslant n+1$ just Erlang-2 inter-observation times suffice (leading to a grand total of $2n^2+n$ moment equations). In that case, define
    \begin{align*}
    \psi_{i,j,k}^{[2]}(\bs\theta)&:=M_{j,k}M_{i,k+2}-\mathbb{E}_{\bs\theta}[M_j(0)M_i(E_{\beta,2})]\ ,
\end{align*}
for $i,j\in\{1,\ldots,n\}$. The system of equations providing the moment estimator $\hat{\bs\theta}_m$ should now also include, besides the $n^2+n$ moment equations \eqref{ME1}--\eqref{ME2} that we already had, additional equations of the type 
\begin{align*}
      \frac{1}{m-2}\sum_{k=1}^{m-2}\psi^{[2]}_{i,j,k}(\hat{\bs \theta}_m) = 0,& \quad i,j\in\{1,\ldots,n \}  ;
\end{align*}
out of these $n^2$ equations the required $(d-1)n$ equations can be picked. In case $d>n+1$, one has to add sufficiently many cross-moments pertaining to higher order Erlang inter-observation times.

To demonstrate how the parametric estimation of the service-time distributions works, we discuss the example in which $G_i$ is distributed exponentially with mean $\mu_i^{-1}.$ It means that $\varrho_i=\lambda_i^{\rm eff}/\mu_i$, and, defining by ${\bs\mu}(\beta)$ the $n$-dimensional vector whose $i$-th entry is $\mu_i/(\mu_i+\beta)$,
\begin{align*}{\boldsymbol P}(\beta) &=  \left(I-{\rm diag}\left\{{\boldsymbol{\mu}(\beta)}\right\}\,Q\right)^{-1} (I-{\rm diag}\{\boldsymbol{\mu}(\beta)\}),\\
{\boldsymbol P}^{\rm (res)}(\beta) &=  {\rm diag}\left\{{\boldsymbol{\mu}(\beta)}\right\}\,Q\,{\boldsymbol P}(\beta)+ (I-{\rm diag}\{\boldsymbol{\mu}(\beta)\});
\end{align*}
here it has been used that for exponentially distributed service times we have that ${\mathscr G}_i(\beta)={\mathscr G}^{\rm (res)}_i(\beta).$
By Proposition \ref{P1}, we can now write $\E[{\bs M}(0){\bs M}(T_\beta)^\top]$ in terms of ${\bs\theta}.$ 
%Along the same lines, $\E[{\bs }(0){\bs N}M(E_{\beta,2})^\top]$ can be expressed in terms of ${\bs\theta}.https://www.overleaf.com/project/661e53ece3973dbbf8669c07#$ 

\liron{
\begin{remark}\label{rem:consistent_G}\em 
   Establishing consistency in this setting requires a somewhat less direct line of argument, since the estimation equation no longer yields an explicit solution. A natural approach is to rely on the theory of Z-estimators (see \cite[Ch.~5]{V2000}). Suppose we impose the requirement that there are $M^-$ and $M^+$ such that, for any $i\in\{1,\ldots,n\}$,
\begin{equation}\label{boundM}0< M^-\leqslant \mu_i\leqslant M^+<\infty.\end{equation}
It would then remain to establish identifiability, i.e. an analogue of Proposition~\ref{prop:identify}, together with suitable uniform regularity of the estimation equations, such as Lipschitz continuity over the parameter space.. Establishing these properties likely requires stronger assumptions. We leave their precise formulation and verification as theoretical open problems for future research and focus here on the practical implementation of the method in this more general setting. $\hfill\Diamond$
\end{remark}

%One could attempt to establish a consistency proof in much the same way as before. Suppose we impose the requirement that there are $M^-$ and $M^+$ such that, for any $i\in\{1,\ldots,n\}$,
%\begin{equation}\label{boundM}0< M^-\leqslant \mu_i\leqslant M^+<\infty.\end{equation}
%Then it is straightforward to verify that under the assumption \eqref{boundM} the analogue of Lemma~\ref{lem:lip} (i.e., Lipschitz continuity)\footnote{LR: remove this reference? What do we want to say here now?} still holds with the new parameter vector $\bs\theta$. However, deriving a counterpart to Proposition~\ref{prop:identify} may prove to be more challenging.
}

A similar approach can be followed as long as the $G_i$ follow a specified distribution that is characterized through a single parameter. Notably, it is not needed that all $G_i$ stem from the same parametric class; one could for instance have that some service times are exponentially distributed, while others are Erlang-2 distributed. In case some of the $G_i$ have distributions with more than one parameter, one needs to add the required number of moment equations in the way discussed above.

\subsection{Model-free approach} \label{ssec:mfree} Interestingly, in case one does not know the parametric form of the service-time distributions, one can still perform an inference procedure. In this case the vector ${\bs\theta}$ consists, besides $Q$ and ${\bs\lambda}$, also of
\[({\mathbb E}[G_1],\ldots,{\mathbb E}[G_n]),\:\:({\mathscr G}_1(\beta),\ldots,{\mathscr G}_n(\beta)),\:\:(\dot{\mathscr G}_1(\beta),\ldots,\dot{\mathscr G}_n(\beta)),\]
i.e., ${\bs\theta}$ has become a vector of dimension $n^2+3n$. Observe that, as a consequence, for $n\geqslant 2$ our $2n^2+n$ moment equations (deriving from the $n$ stationary moments, the $n^2$ cross-moments pertaining to exponential inter-observation times, and the $n^2$ cross-moments pertaining to Erlang-2 inter-observations times) suffice.

%It is directly seen that $\bs\varrho$, $\E[{\bs M}(0){\bs M}(T_\beta)^\top]$, and  $\E[{\bs M}(0){\bs M}(E_{\beta,2})^\top]$ can be expressed in terms of the entries of $\bs\theta$. Assuming that,  for any $i\in\{1,\ldots,n\}$,  ${\mathbb E}[G_i]$,  ${\mathscr G}_i(\beta)$, and $\dot{\mathscr G}_i(\beta)$ are Lipschitz in ${\bs\theta}\in\Theta$ and $0<g^-\leqslant{\mathbb E}[G_i]\leqslant g^+<\infty$, one again obtains consistency; the requirement ${\mathbb E}[G_i]\geqslant g^->0$ is needed to make sure that 
%\[\mathscr{G}_i^{\rm (res)}(\beta) = \frac{1-\mathscr{G}_i(\beta)}{\beta\, \E[G_i]}\quad\mbox{and}\quad \dot{\mathscr{G}}_i^{\rm (res)}(\beta)= -\frac{1}{{\mathbb E}[G_i]}\frac{1-{\mathscr G}_i(\beta)+\beta \,\dot{\mathscr G}_i(\beta)}{\beta^2}\]
%are Lipschitz for $\bs\theta\in\Theta$.
%\footnote{LR: I removed the entire paragraph here because a discussion will essentially be identical to the previous subsection.}

Evidently, this approach is particularly useful in case one does not know the service-time distribution's parametric form. The price to be paid is that one needs to estimate as many as three parameters per service time.

\section{Unknown observation probabilities}
\label{sec:obsprob}

\michel{So far we have assumed that the full network population vector is observed. In this section, we consider the more general setting of {\it censored observations}, in which at any timestep only a subset of customers is observed. Besides the parameters of the infinite-server queueing network, this setting introduces additional parameters describing the underlying observation mechanism.
Our main objective in this section is to derive the corresponding moment expressions for the {\it observed} network population. These identities provide the basis for extending the estimation framework developed in the previous sections to jointly estimate both the network parameters and the observation probabilities. While implementing the extension is straightforward and demonstrates good performance, as will be shown in the sequel, rigorous verification of consistency requires additional identifiability conditions and a separate consistency analysis, which lie beyond the scope of the present paper.}

Let $p_j$ be the probability that a customer residing at station $j$ is actually observed. We assume that each customer is observed independently, i.e., if there happen to be $m_j\in{\mathbb N}_0$ customers at station $j$, the number of observed customers is binomially distributed with parameters $m_j$ and $p_j$. In addition, we assume that this observation process happens `independently in time' (i.e., at every new observation time, the `binomial selection' of which customers are observed occurs independently from the selection at previous observation times) as well as `independently in space' (i.e., at every station, the `binomial selection' occurs independently from the selection at other stations).

\michel{By the tower property, the mean number of observed customers at station $i$ at time $t$, denoted by $N_i(t)$, is, for any $t\geq 0$, given by}
\[
\E[N_i(t)]
=\E\big[\E[N_i(t)\mid {\boldsymbol M}(t)]\big]
=\E[p_i M_i(t)]
=p_i\varrho_i
=p_i\lambda_i^{\rm eff}\E[G_i].
\]
Also, again by the tower property, at any $t\geqslant 0$,
\begin{align*}
\E[N_j(t)N_i(t+T_\beta)]&=\E[N_j(0)N_i(T_\beta)]  = \E\big[\E[N_j(0)N_i(T_\beta)\,|\,{\boldsymbol M}(0),{\boldsymbol M}(T_\beta)]\big]\\
    &= \E\big[p_j p_iM_j(0)M_i(T_\beta)\big] = p_j p_i\,\E\big[M_j(0)M_i(T_\beta)\big] ,
\end{align*}
with $\E\big[M_j(0)M_i(T_\beta)\big]$ as given in Proposition \ref{P1}. We have shown the following statement.

\begin{lemma}
    \label{L3}
  For any $\beta>0$,
  \[\E[{\bs N}(0){\bs N}(T_\beta)^\top]={\rm diag}\{{\bs p}\}\,\E[{\bs M}(0){\bs M}(T_\beta)^\top]\,{\rm diag}\{{\bs p}\},\]
  with $\E[{\bs M}(0){\bs M}(T_\beta)^\top]$ being given by Proposition~\ref{P1}.
\end{lemma}

% Usning this we can write an expression for $\E(M_i(t_0 + T_\beta)|M_j(t_0) = m_j)$,
% \begin{align*}
%     \E(M_i(t_0 + T_\beta)|M_j(t_0) = m_j) &= m_j\, P_{ji}^{\rm (res)}(\beta) + \sum_{j=1}^n \Gamma_{ji}(\beta) + \sum_{k \neq j} \sum_{m_k = 0}^\infty m_k \,\mathbb{P}(M_k(t_0) = m_k) \,P_{ki}^{\rm (res)}(\beta)\\
%     &=  m_j\, P_{ji}^{\rm (res)}(\beta) + \sum_{j=1}^n \Gamma_{ji}(\beta) + \sum_{k \neq j} \lambda_k^{eff} \,t_0 \,P_{ki}^{\rm (res)}(\beta).
% \end{align*}
% Now, we can use the expression of the conditional expectation to find the cross-moments
% \[\E(M_i(T_\beta)M_j(t_0)) = \sum_{m_j = 0}^\infty m_j \, \mathbb{P}(M_j(t_0) = m_j)\,\E(M_i(T_\beta)|M_j(t_0) = m_j).\] 
Because the number of unknown parameters exceeds $n^2+n$ (even in the case the service-time distributions are known), we need additional moment equations. We could again utilize the moment equations deriving from $\E[{\bs N}(0){\bs N}(E_{\beta,2})^\top]$, where it is noted that this covariance matrix can be found from $\E[{\bs N}(0){\bs N}(T_\beta)^\top]$ in the same way as $\E[{\bs M}(0){\bs M}(E_{\beta,2})^\top]$ was found from $\E[{\bs M}(0){\bs M}(T_\beta)^\top]$.

\section{Numerical experiments}\label{sec:exp}
In this section we assess the performance of our estimator through a series of experiments. Before doing so, we discuss the moment equations that are to be solved in Subsection \ref{ssec:momeq}, and provide a specific algorithm for the case that $\bs p=\bs 1$ in Subsection \ref{ssec:p1}. The actual experiments are presented in Subsection \ref{ssec:exp}.

\subsection{Solving the moment equations}\label{ssec:momeq}
To evaluate our estimator, it is crucial to ensure that the number of equations matches the number of unknown parameters. \michel{Two procedures have been used.}

\michel{\begin{itemize}
\item[(P1)]
In case the unknown parameters are ${\bs\lambda}$ and $Q$, as in Section \ref{sec:method}, we can base our computations on the moment equations related to the `snapshot means' $\E[{\bs M}(0)]$ and lag-one covariances $\E[{\bs M}(0){\bs M}(T_\beta)^\top]$. This leads to $n^2+n$ equations that define the $n^2+n$ parameters, where the estimator follows directly from Proposition~\ref{prop:identify}. 
As described in Subsection~\ref{ssec:est}, the truncation~\eqref{eq:lambda_truncated}--\eqref{eq:q:truncated} is applied to make sure that the resulting estimator lies in $\Theta$.
\item[(P2)]
In case there are more parameters, we have to add sufficiently many moment equations related to lag-two covariances. For instance, in the case that 
the service-time distributions \( G_i \) are characterized by a single parameter $\eta_i$ and the observation probabilities $p_i$ are unknown, there are \( n^2 + 2n \) unknown parameters: recalling that in this framework the diagonal elements of $Q$ are assumed $0$, there are $n^2-n$ parameters related to $Q$, and $n$ parameters related to each of the vectors $\bs\lambda$, ${\bs\eta}$ and $\bs p$.  This requires us to use the $n$ `snapshot means' $\E[{\bs N}(0)]$, the $n^2$ lag-one covariances $\E[{\bs N}(0){\bs N}(T_\beta)^\top]$ and in addition $n$ lag-two covariances $\E[{\bs N}(0){\bs N}(E_{\beta,2})^\top]$ (where we picked the cross moments that appear on the diagonal of this matrix). The parameter estimates are obtained by equating these analytical expressions of the moments (which are functions of the unknown parameters) to their empirical counterparts. 

\noindent
Observe that simply solving the moment equations does not guarantee that the solution lies within the parameter space $\Theta$; in principle, it may lead to, e.g., negative arrival rates, negative probabilities, or row sums of the routing matrix exceeding 1. 
In our numerical experiments, we therefore adopted a pragmatic approach: minimizing the sum of squared differences between the left-hand and right-hand sides of the moment equations, with the optimization carried out over ${\bs\theta} \in \Theta$.
\end{itemize}}
\subsection{Efficient procedure if all customers are observed}\label{ssec:p1}

Supposing that ${\bs p = \bs 1}$ and again considering the situation that each service-time distribution is characterized by a single parameter, we have $n^2+n$ unknown parameters. The parameters (i.e., the entries of $Q$, $\bs{\eta}$ and $\bs{\lambda}$) can be estimated using the moment equations in Equation~\eqref{eq:loads2} and Proposition~\ref{P1}, replacing $\E[{\bs M}(0)]$ and $\E[{\bs M}(0){\bs M}(T_\beta)^\top]$ with their empirical counterparts. Specifically, define the empirical first-moment vector, for $i\in\{1,\ldots,n\}$,
$\bs{\psi}^{[0]}$ by \[\bs{\psi}^{[0]}_i = \frac{1}{m}\sum_{k=1}^m\psi^{[0]}_{i,k}(\bs\theta),\] and the empirical second-moment matrix $\bs{\psi}^{[1]}$ by, for $i,j\in\{1,\ldots,n\}$,  \[\bs{\psi}^{[1]}_{ij} =\frac{1}{m-1} \sum_{k=1}^{m-1}\psi^{[1]}_{i,j,k}(\bs\theta).\] Let $\bs G$ denote the vector of expected service times, i.e., $\bs{G}_i = \E[G_i]$, and define
\begin{align*}
    t_1:&=  \bs{\psi}^{[0]}\, (\bs{\psi}^{[0]})^\top +{\rm diag}\{\bs{\psi}^{[0]}\}, \\
    t_2:&={\beta}^{-1}\,\bs{\psi}^{[0]}\, (\bs{\psi}^{[0]})^\top\, \text{diag} \{{\bs G }\}^{-1} \, (I-Q)  , \\
    t_3: &= (I - \text{diag}\{{\boldsymbol{\mathscr G}}(\beta)\} Q)^{-1} \big( I - \text{diag}\{{\boldsymbol{\mathscr G}}(\beta)\} \big) \,.
\end{align*}
The system of equations to solve is
\begin{align*}
    & \bs \psi^{[0]} =  {\rm diag}\{\bs G\} \, (I-Q^\top)^{-1}{\boldsymbol \lambda}\\
    & \bs \psi^{[1]} = t_1\,\big(I- {\rm diag}\{{\boldsymbol{\mathscr G}^{\rm (res)}}(\beta)\}\big) + \big(t_1 \, {\rm diag}\{{\boldsymbol{\mathscr G}^{\rm (res)}}(\beta)\}\,Q + t_2 \big)\,t_3 \,.
\end{align*}
Here we should realize that ${\bs G}$, ${\mathscr G}(\beta)$ and ${\mathscr G}^{\rm (res)}(\beta)$ can all be expressed in terms of the $n$-dimensional unknown parameter vector $\bs\eta$.

Note that in the above moment equations we have used $\bs \psi^{[0]}$ in place of the explicit expression for $\bs \varrho$ in the cross-moment equations, so as to simplify the computations.
Furthermore, observe that the cross-moment equations can be written solely in terms of $Q$ and $\bs \eta$; indeed, we can eliminate $\bs\lambda$ due to the relation
\[
    {\bs \lambda}^\top =  \big({\bs \lambda}^{\rm eff}\big)^\top \, (I-Q) = {\bs \varrho}^\top \, \text{diag} \{{\bs G }\}^{-1} \, (I-Q) \,.
\]
After having eliminated the $\bs\lambda$, we have arrived at a system of $n^2$ equations in the $n^2$ unknowns (i.e., the parameters relating to $Q$ and ${\bs\eta}$). Hence we can  solve for $Q$ and $\bs \eta$ using the cross-moment equations, and then subsequently identify $\bs \lambda$ from the single snapshot equations using the estimated $\hat{Q}$ and $\hat{\bs \eta}$.
We finally remark that the right-hand side of the cross-moment equations is quite intricate and can be computationally demanding. Our numerical experiments indicate that the following alternative formulation, obtained by an elementary reordering, often leads to a significant reduction in computation time:
\[   
\big(t_1\, \text{diag}\{{\boldsymbol{\mathscr G}^{\rm (res)}}(\beta)\} Q + t_2\big)^{-1} \, \big( \bs \psi^{[1]} - t_1 \big(I - \text{diag}\{{\boldsymbol{\mathscr G}^{\rm (res)}}(\beta)\}\big) \text{diag} \{{\bs p }\}   \big) - t_3= 0 \,.
\]

\subsection{Experiments}\label{ssec:exp}
In this subsection, we present a series of experiments that demonstrate the effectiveness of our estimation procedure. These experiments address various aspects of the estimator: in Experiment 1, we assess the impact of the network structure; in Experiment 2, we examine how performance accuracy varies with the number of observations; in Experiment 3, we show that the procedure remains effective even for large-scale networks; and in Experiment 4, we illustrate how the method can be applied when the service-time distributions have no known parametric form.

\medskip

\begin{figure}[htbp]
    \centering
    \begin{subfigure}[b]{0.23\textwidth}
        \includegraphics[width=\textwidth]{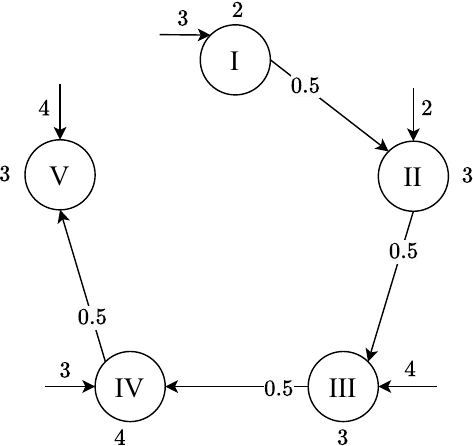}
        \caption*{(a)~{\it Line}}
    \end{subfigure}
     \hfill % Add horizontal fill for spacing
    \begin{subfigure}[b]{0.23\textwidth}
        \includegraphics[width=\textwidth]{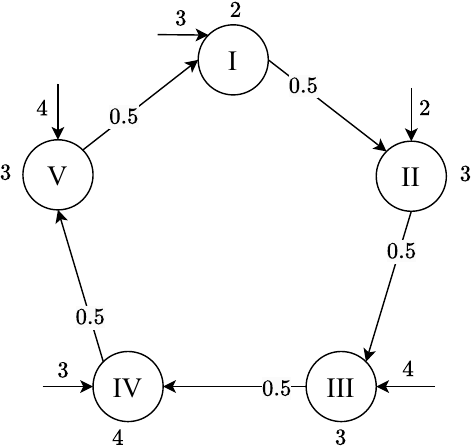}
        \caption*{(b)~{\it Circle}}
        \label{fig:1b}
    \end{subfigure}
    % \vspace{0.5cm} % Vertical spacing between rows
    \hfill
    \begin{subfigure}[b]{0.23\textwidth}
        \includegraphics[width=\textwidth]{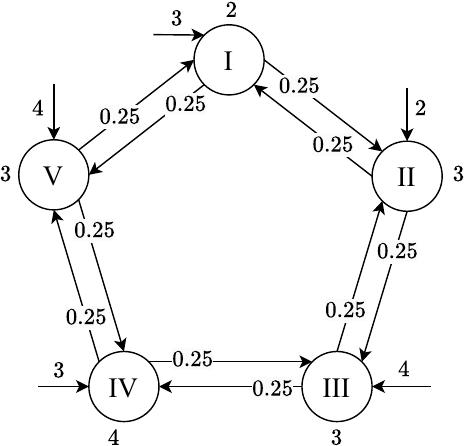}
        \caption*{(c)~{\it Symmetric Circle}}
        \label{fig:sub3}
    \end{subfigure}
    \hfill
    \begin{subfigure}[b]{0.23\textwidth}
        \includegraphics[width=\textwidth]
        {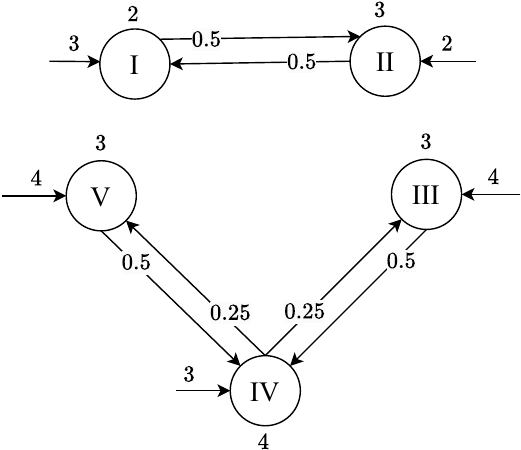}
        \caption*{(d)~{\it 2+3 Cliques}}
    \end{subfigure}
    \caption{Four topologies considered in Experiment 1.}
    \label{fig:four_figures}
\end{figure}

{\it Experiment~1.} In the first experiment, we consider an instance in which we  wish to estimate the vector of external arrival rates $\bs\lambda$ and the transition rate matrix $Q$; the service times are exponentially distributed with {known} rates $\bs\mu$, and the observation probabilities $\bs p$ are known to equal 1. It means that we are in a setting in which Proposition \ref{prop:identify} can be applied \michel{(in combination with the truncation~\eqref{eq:lambda_truncated}--\eqref{eq:q:truncated}; see (P1) in Subsection~\ref{ssec:momeq})}. 
In our experiment, we considered a system of \( n = 5 \) stations, where the queue-length vector \( \bs{M}(\cdot) \) was observed at Poisson-distributed time intervals with rate \( \beta = 5 \). The system parameters were defined as follows: the arrival rates are given by the vector \( \bs{\lambda} = (3, 2, 4, 3, 4)^\top \), and the service rates by \( \bs{\mu} = (2, 3, 3, 4, 3)^\top \). We study four distinct network topologies, each characterized by a different routing matrix \({Q} \):  
\begin{itemize}
    \item[(a)] {\it Line}, with \( q_{i,i+1} = 0.5 \) for \( i = 1, \dots, 4 \), whereas \( q_{ij} = 0 \) otherwise. 
    \item[(b)] {\it Circle}, with \( q_{i,i+1} = 0.5 \) for \( i = 1, \dots, 4 \) and \( q_{5,1} = 0.5 \), whereas \( q_{ij} = 0 \) otherwise.  
    \item[(c)] {\it Symmetric Circle}, with \( q_{i,i+1} = q_{i+1,i} = {0.25} \) for \( i = 1, \dots, 4 \), \( q_{1,5} = q_{5,1} = {0.25} \), whereas \( q_{ij} = 0 \) otherwise.   
    \item[(d)] {\it 3 + 2 Cliques},  with \( q_{1,2} = q_{2,1} = 0.5 \), \( q_{3,4} = q_{5,4} = 0.5 \),  \( q_{4,3} = q_{4,5} = 0.25 \), whereas \( q_{ij} = 0 \) otherwise.  
\end{itemize}
The four topologies are schematically illustrated in Figure \ref{fig:four_figures}. The numbers on the arrows between stations represent routing probabilities, those on the arrows entering the network indicate external arrival rates, and the numbers next to the stations denote service rates.

% Hence, the underlying chain structure is that of a tandem system\footnote{LR: Add a diagram of the true network? It might be a good idea to add diagrams of all the true networks and number them (Network 1,2,3... or Example 1,2,3...). Then in the estimation results you can reference the network you are dealing with.}, and 
The first objective of this experiment is to verify, for each of the four topologies, whether our estimation procedure can identify the structure of the underlying network. In the second place, we wish to study the precision of the corresponding parameter estimates. 
To this end, we performed for each of the four topologies \hritika{$R=100$} runs, each comprising \hritika{$m=2\,500\,000$} observation times.
We thus compute sample means and sample variances based on \hritika{$100$} estimates of \(\hat{Q}\) and \(\hat{\lambda}\). For each of the four instances, these are summarized below.

\begin{itemize}
\item[(a)] In the line network the sample mean and sample variance, based on the $R$ values of \(\hat{\bs\lambda}\), are  given by
  \begin{align*}{\bs m}_{\bs\lambda} &= 
  (2.9572, 1.9643, 3.9571, 2.9676, 3.9581) \\ {{\bs v}_{\bs\lambda}}&=(0.0006, 0.0010, 0.0013, 0.0010, 0.0016),
 \end{align*}
  respectively. The sample mean and variance of \(\hat{Q}\), based on the $R$ runs, are
  \begin{align*}{\bs m}_{Q}&= 
  \begin{pmatrix*}[r]
  0.0035 & 0.5011 & 0.0022 & 0.0024 & 0.0031 \\
  0.0026 & 0.0020 & 0.5008 & 0.0016 & 0.0025 \\
  0.0015 & 0.0012 & 0.0022 & 0.5000 & 0.0018 \\
  0.0009 & 0.0014 & 0.0018 & 0.0017 & 0.4999 \\
  0.0013 & 0.0015 & 0.0015 & 0.0013 & 0.0019
  \end{pmatrix*},\\
  {\bs v}_{Q}&= 10^{-4} \cdot
  \begin{pmatrix*}[r]
  0.2217 & 0.2683 & 0.1154 & 0.1008 & 0.1720 \\
  0.1080 & 0.0942 & 0.2949 & 0.0803 & 0.1181 \\
  0.0432 & 0.0372 & 0.0811 & 0.1623 & 0.0782 \\
  0.0220 & 0.0405 & 0.0576 & 0.0617 & 0.1890 \\
  0.0318 & 0.0252 & 0.0601 & 0.0411 & 0.1076
  \end{pmatrix*},\end{align*}
  respectively.
  The experiments thus reveal that the estimation algorithm typically picks up the network structure, and that the standard deviations of all $30$ parameters (being the square roots of the sample variances) of the estimator are low relative to the estimated value. 
  
\item[(b)] In the circle network, the sample mean and sample variance of \(\hat{\bs\lambda}\) are
\begin{align*}{\bs m}_{\bs\lambda} &= 
  (2.9444, 1.9569, 3.9514, 2.9587, 3.9611) \\ {\bs v}_{\bs\lambda}&=(0.0021, 0.0012, 0.0016, 0.0016, 0.0015),
 \end{align*}
  respectively. The sample mean and variance of \(\hat{Q}\) are
  \begin{align*}{\bs m}_{Q}&= 
  \begin{pmatrix*}[r]
  0.0026 & 0.5005 & 0.0021 & 0.0023 & 0.0015 \\
  0.0025 & 0.0023 & 0.5004 & 0.0014 & 0.0017 \\
  0.0022 & 0.0016 & 0.0019 & 0.5001 & 0.0016 \\
  0.0020 & 0.0013 & 0.0014 & 0.0017 & 0.4992 \\
  0.4998 & 0.0012 & 0.0015 & 0.0011 & 0.0020
  \end{pmatrix*},\\
  {\bs v}_{Q}&= 10^{-4} \cdot
  \begin{pmatrix*}[r]
   0.1517 & 0.1824 & 0.0755 & 0.0804 & 0.0608 \\
  0.1105 & 0.1153 & 0.2343 & 0.0706 & 0.0623 \\
  0.0929 & 0.0533 & 0.0801 & 0.1867 & 0.0574 \\
  0.0651 & 0.0510 & 0.0468 & 0.0507 & 0.1585 \\
  0.1605 & 0.0372 & 0.0642 & 0.0333 & 0.0915
  \end{pmatrix*},\end{align*}
  respectively. \hritika{Again, the network structure is picked up well.}
  % Also observe that this experiment indicates that our estimation procedure can effectively distinguish this topology, in which the clients are routed in a clockwise manner, from its counterpart in which the clients are routed in a counterclockwise manner (with the same parameters), despite the two systems having the same equilibrium distribution. 

\item[(c)] In the symmetric circle network, the sample mean and sample variance of \(\hat{\bs\lambda}\) are
\begin{align*}{\bs m}_{\bs\lambda} &= 
  (2.9554, 1.9679, 3.9590, 2.9694, 3.9524)\\ {\bs v}_{\bs\lambda}&=(0.0028, 0.0024, 0.0022, 0.0018, 0.0028),
 \end{align*}
 respectively.
The sample mean and variance of \(\hat{Q}\) are
 \begin{align*}{\bs m}_{Q}&= 
  \begin{pmatrix*}[r]
 0.0031 & 0.2496 & 0.0019 & 0.0018 & 0.2507 \\
  0.2500 & 0.0026 & 0.2499 & 0.0018 & 0.0019 \\
  0.0018 & 0.2499 & 0.0027 & 0.2499 & 0.0017 \\
  0.0016 & 0.0014 & 0.2501 & 0.0013 & 0.2508 \\
  0.2506 & 0.0016 & 0.0018 & 0.2503 & 0.0021
  \end{pmatrix*},\\
  {\bs v}_{Q}&= 10^{-4} \cdot
  \begin{pmatrix*}[r]
  0.2021 & 0.2261 & 0.0863 & 0.0613 & 0.3118 \\
  0.3000 & 0.1660 & 0.1999 & 0.0590 & 0.0996 \\
  0.0775 & 0.1268 & 0.1171 & 0.1327 & 0.0519 \\
  0.0509 & 0.0394 & 0.1626 & 0.0389 & 0.1832 \\
  0.2102 & 0.0450 & 0.0674 & 0.1719 & 0.1056
  \end{pmatrix*},\end{align*}
  respectively. As was the case in the previous instances, we observe that the procedure is capable of identifying the system's structure. \hritika{We also note that, in this case, the sample variances of \(\hat{\bs\lambda}\) are slightly higher than those for the line network.}

  \item[(d)] In the 3 + 2 cliques network, the sample mean and sample variance of \(\hat{\bs\lambda}\) are
  \begin{align*}{\bs m}_{\bs\lambda} &= 
  (2.9503, 1.9597, 3.9567, 2.9623, 3.9603)\\ {\bs v}_{\bs\lambda}&=(0.0019, 0.0012, 0.0015, 0.0025, 0.0015),
 \end{align*}
  respectively. The sample mean and variance of \(\hat{Q}\) are
   \begin{align*}{\bs m}_{Q}&= 
  \begin{pmatrix*}[r]
 0.0025 & 0.4996 & 0.0023 & 0.0025 & 0.0019 \\
  0.5000 & 0.0023 & 0.0018 & 0.0019 & 0.0021 \\
  0.0020 & 0.0014 & 0.0020 & 0.4999 & 0.0015 \\
  0.0015 & 0.0012 & 0.2503 & 0.0019 & 0.2497 \\
  0.0016 & 0.0016 & 0.0013 & 0.4998 & 0.0020
  \end{pmatrix*},\\
  {\bs v}_{Q}&= 10^{-4} \cdot
  \begin{pmatrix*}[r]
0.1750 & 0.2339 & 0.1083 & 0.1639 & 0.1004 \\
  0.3393 & 0.1159 & 0.0725 & 0.0803 & 0.0715 \\
  0.0806 & 0.0441 & 0.1090 & 0.2576 & 0.0542 \\
  0.0522 & 0.0283 & 0.1337 & 0.0696 & 0.1109 \\
  0.0488 & 0.0437 & 0.0571 & 0.2000 & 0.0873
  \end{pmatrix*},\end{align*}
  respectively. The cliques can be effectively identified from the estimated matrix $\hat{Q}$.
  \end{itemize}

\begin{figure*}[h!]
\centering
\begin{subfigure}[t]{0.32\textwidth}
    \centering
    \includegraphics[width=1\linewidth]{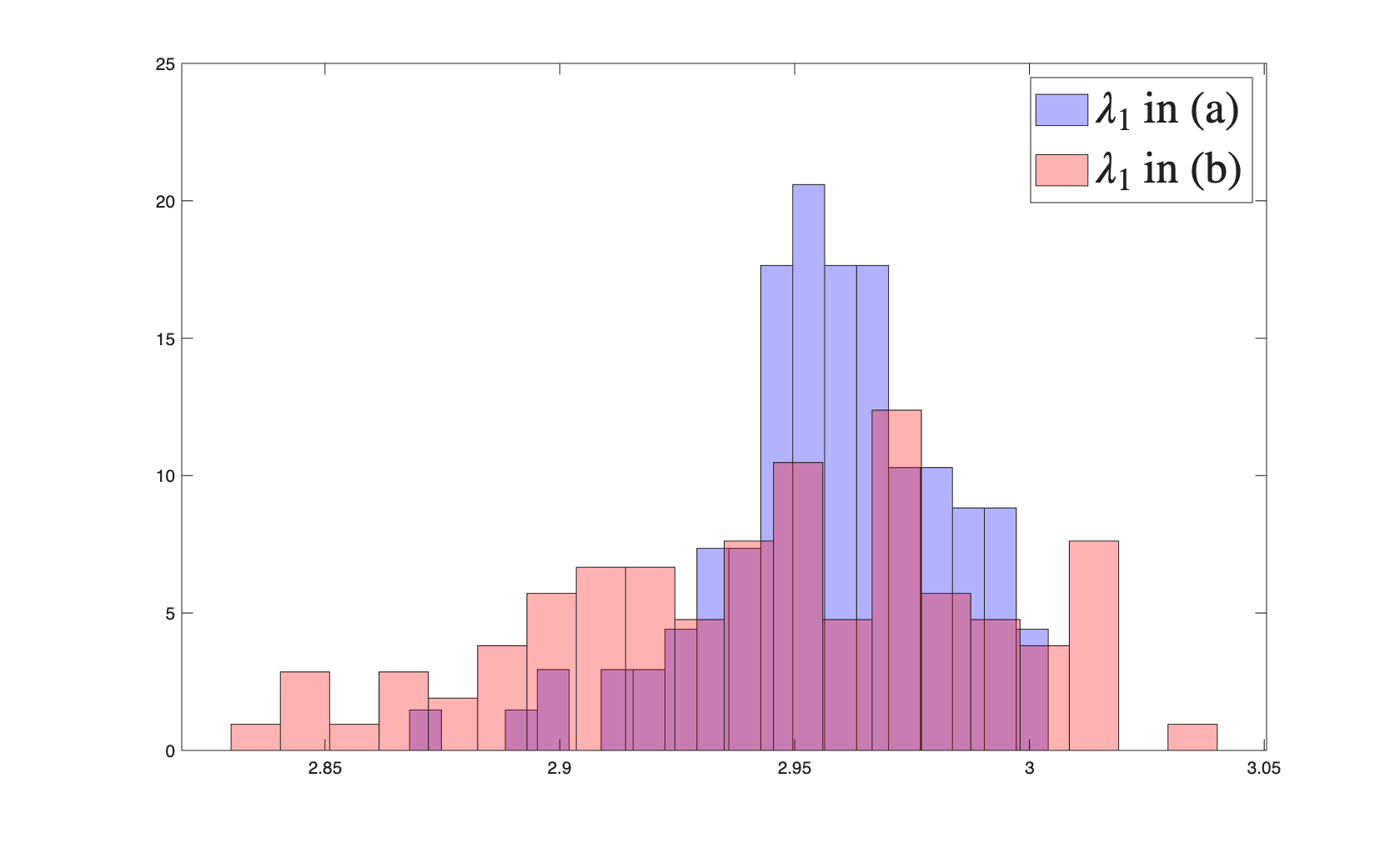}
    \label{fig:lambda1ab}
\end{subfigure}
\hfill
\begin{subfigure}[t]{0.32\textwidth}
    \centering
    \includegraphics[width=1\linewidth]{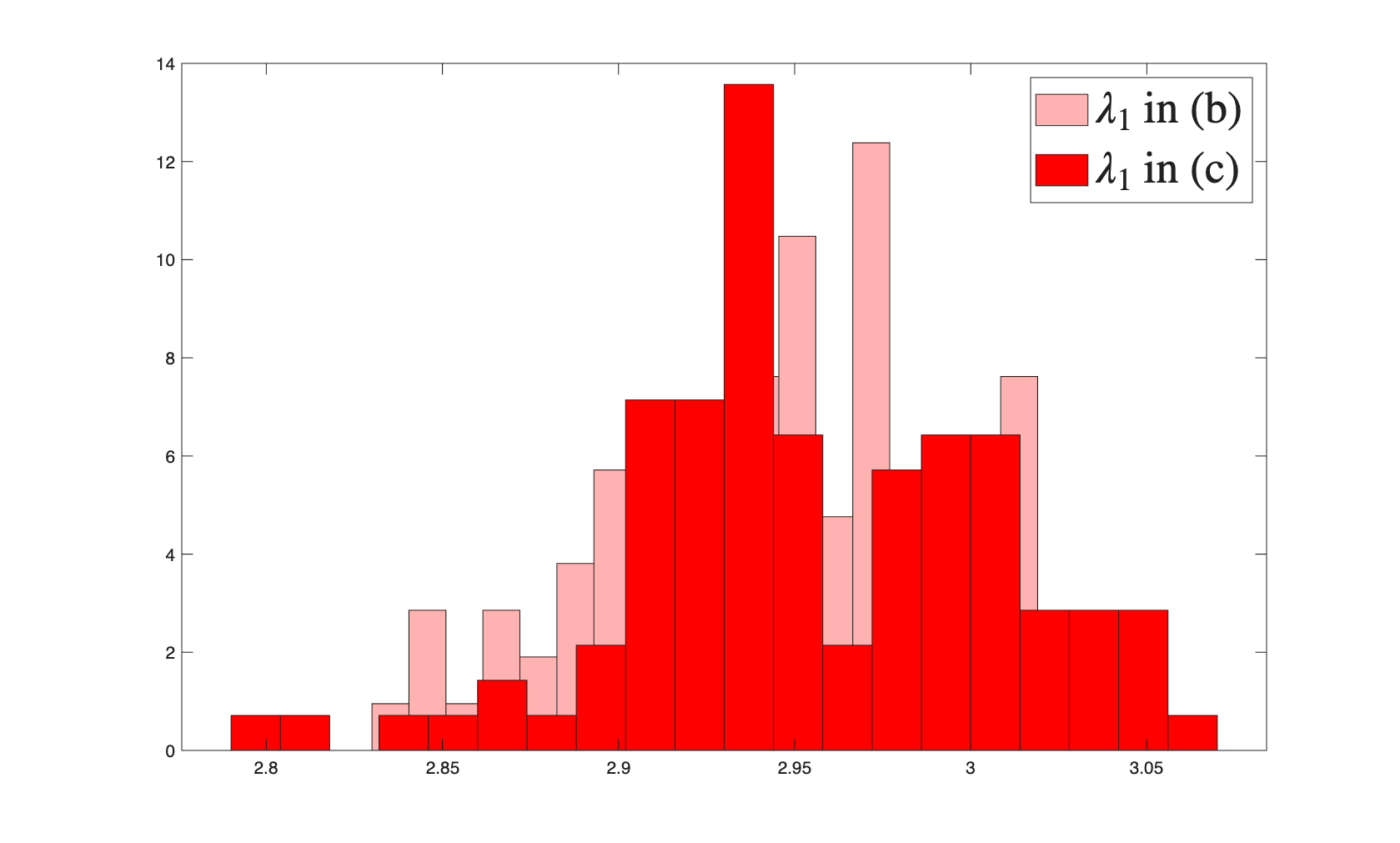}
    \label{fig:lambda1bc}
\end{subfigure}
\hfill
\begin{subfigure}[t]{0.32\textwidth}
    \centering
    \includegraphics[width=1\linewidth]{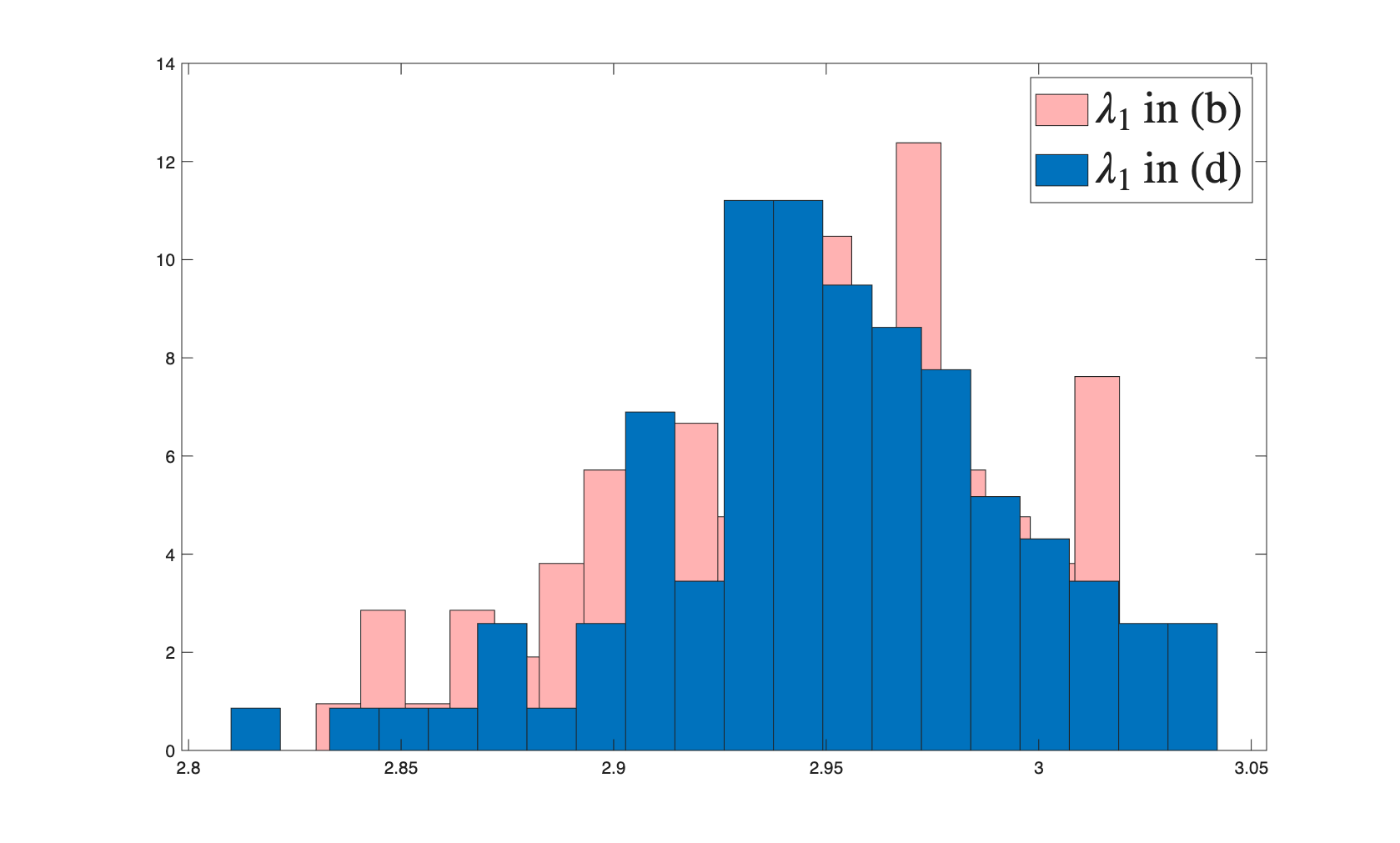}
    \label{fig:lambda1bd}
\end{subfigure}
% \begin{subfigure}[t]{0.48\textwidth}
%     \centering
%     \includegraphics[width=1\linewidth]{q11.jpg}
%     \caption{Histogram for $\hat q_{11}$}
%     \label{fig:q11}
% \end{subfigure}
% \hfill
% \begin{subfigure}[t]{0.48\textwidth}
%     \centering
%     \includegraphics[width=1\linewidth]{q12.jpg}
%     \caption{Histogram for $\hat q_{12}$}
%     \label{fig:q12}
% \end{subfigure}
\caption{\label{Fig1}Histograms pertaining to estimates of $\lambda_1$, based on \hritika{$R=100$} experiments, each with \hritika{$m=2\,500\,000$} observations.}
\end{figure*}

The above findings confirm the desirable property  that, for the four instances considered, the standard deviations of the estimator are low relative to the estimated value. We conclude this experiment with a somewhat more detailed analysis of the distribution of our estimator. 
\hritika{Figure~\ref{Fig1}} displays histograms of the estimates for \( \lambda_1 \); the estimates for \( \lambda_i \) (\( i = 2,\ldots,5 \)) and the entries of $Q$ exhibit similar behavior. The results demonstrate that the distribution of the estimator is, in all cases, centered around the true value. \michel{The experiments show that the variances of \( \hat{\lambda}_i \), \(i=1,\ldots,5\), are of similar order across all four topologies. The estimates of \( \hat{Q} \) likewise exhibit consistently small variances, with only minor differences across the four topologies. This behavior may be explained by the fact that, in all four cases, the same number of entries of the \(Q\)-matrix is estimated; hence, there is no particular advantage to exploiting sparsity in \(Q\).
} 

\begin{remark}\em 
    An important strength of our technique is its ability to distinguish between different network topologies, even when the underlying stochastic systems share the same stationary distribution. For example, consider the following two cases: (i) a circular topology with homogeneous arrival rates ($\lambda_i = \lambda$ for all $i\in\{1,\ldots,n\}$) and homogeneous service rates ($\mu_i = \mu$), with clockwise routing and uniform exit probabilities (say, $q$); and (ii) an otherwise identical system, but with counterclockwise routing. Although these two systems clearly have the same stationary distribution (which is $n$-dimensional Poisson), our estimation procedure can successfully distinguish between the two underlying topologies. \hfill $\Diamond$
\end{remark}

% \textit{Experiment 1.2} Now we present a similar analysis as above for a circular network. The goal is compare the variances and distribution of the estimates for the two network. $\lambda$ and $\mu$ are same as above. $Q$ is such that $q_{i,i+1} = 0.5$ for $i = 1,\dots,4$, $q_{51} = 0.5$ and $q_{ij} = 0$ else. 

\begin{figure*}[t!]
\centering
\begin{subfigure}[t]{0.48\textwidth}
    \centering
        \includegraphics[clip, trim=5.5cm 18.5cm 4.5cm 3cm, width=1\textwidth]{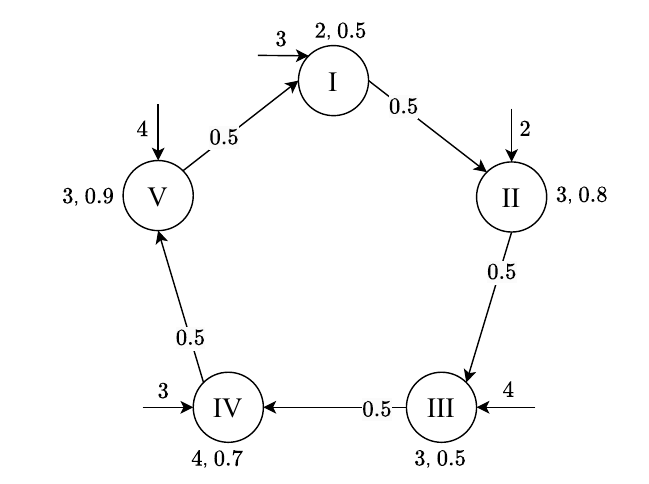}
    \caption*{(a)~True network.}
    \label{fig:example_true}
\end{subfigure}
\hfill
\begin{subfigure}[t]{0.48\textwidth}
    \centering
        \includegraphics[clip, trim=5.5cm 18.5cm 4.5cm 3cm, width=1\textwidth]{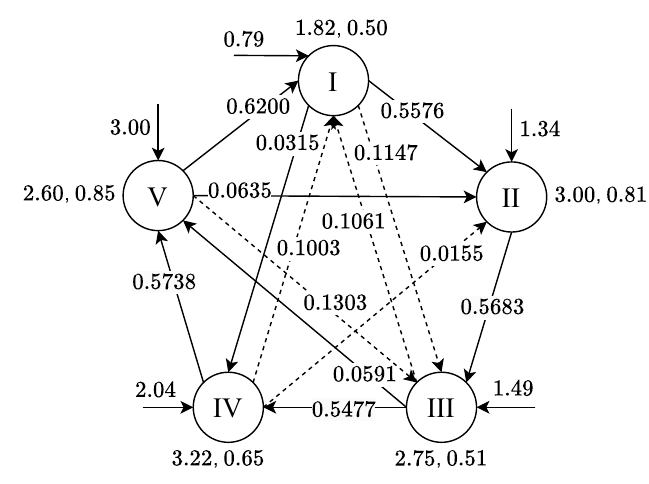}
    \caption*{(b)~Estimated network, $m=10^5$.}
    \label{fig:example_hundred}
\end{subfigure}

\begin{subfigure}[t]{0.48\textwidth}
    \centering
        \includegraphics[clip, trim=5.5cm 18.5cm 4.5cm 3cm, width=1\textwidth]{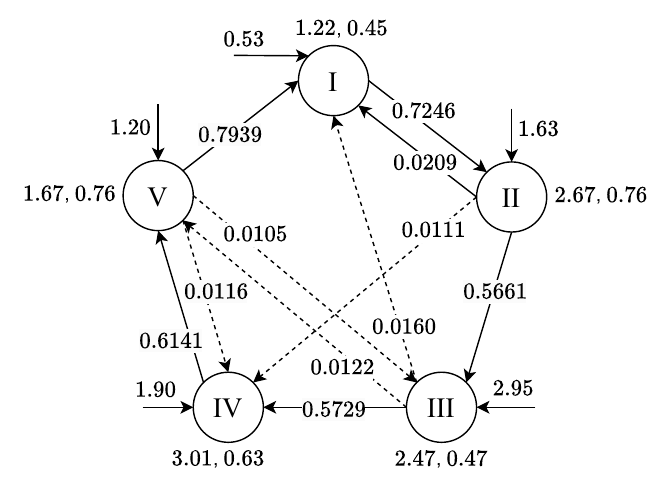}
    \caption*{(c)~Estimated network, $m=5\cdot10^5$.}
    \label{fig:example_fivehundred}
\end{subfigure}
\hfill
\begin{subfigure}[t]{0.48\textwidth}
    \centering
        \includegraphics[clip, trim=5.5cm 18.5cm 4.5cm 3cm, width=1\textwidth]{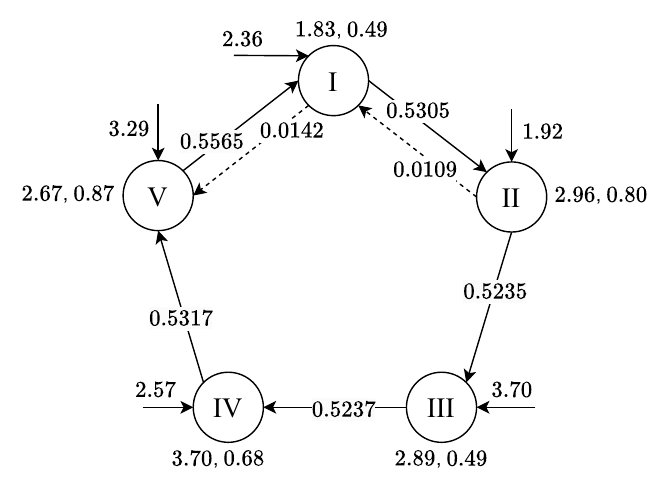}
    \caption*{(d)~Estimated network, $m=10^6$.}
    \label{fig:example_onemillion}
\end{subfigure}

\begin{subfigure}[t]{0.48\textwidth}
    \centering
        \includegraphics[clip, trim=5.5cm 18.5cm 4.5cm 3cm, width=1\textwidth]{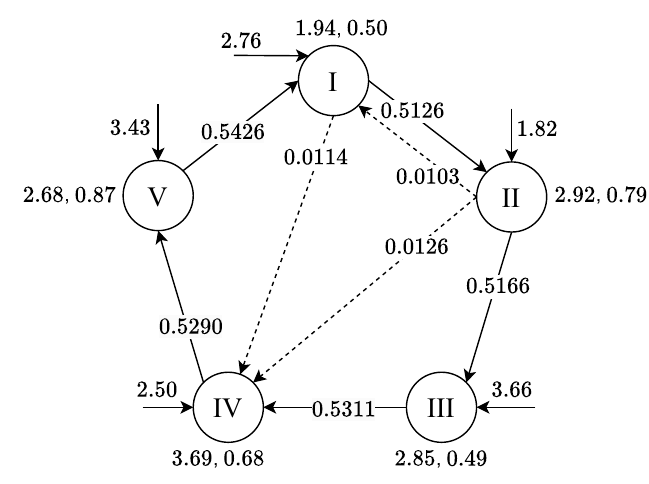}
    \caption*{(e)~Estimated network, $m=2\cdot 10^6$.}
    \label{fig:example_twomillion}
\end{subfigure}
\hfill
\begin{subfigure}[t]{0.48\textwidth}
    \centering
        \includegraphics[clip, trim=5.5cm 18.5cm 4.5cm 3cm, width=1\textwidth]{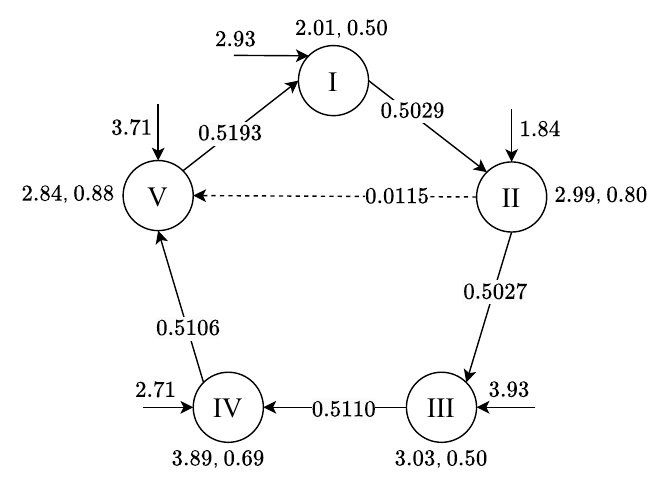}
    \caption*{(f)~Estimated network, $m=3\cdot 10^6$.}
    \label{fig:example_threemillion}
\end{subfigure}
\caption{An example with five stations and a cyclic network.\label{Fig2}}
\end{figure*}

{\it Experiment 2.} The main goal of this experiment is to analyze the impact of the number of observations on the accuracy of the estimates. We consider a network consisting of five stations and the routing matrix is cyclic, as depicted in Panel~({a}) of Figure \ref{Fig2}. The service times are exponentially distributed with unknown parameter vector ${\bs \mu}$. In this experiment we are estimating the $30$ parameters corresponding to ${\bs \lambda}$, $Q$, ${\bs \mu}$, and the observation probability vector ${\bs p}$ (being in a setup in which the diagonal elements of $Q$ are equal to $0$). In Panel~({a}) of Figure~\ref{Fig2}, the two numbers next to station~$i$, for $i={\rm I},\ldots,{\rm V}$, represent $\mu_i$ and $p_i$, respectively. The routing probabilities $q_{ij}$ are indicated by the numbers on the arrows between stations, while the external arrival rates $\lambda_i$ are shown on the arrows pointing into the stations.

Panels~({b})-({f}) of Figure~\ref{Fig2} provide the estimates for different numbers of observations $m$, with the sampling being performed at rate $\beta = 5$. For the estimated network, an arrow has been omitted if the value of the corresponding routing probability is below 0.01, and is dotted if it is below 0.02. 
%\footnote{\scriptsize \tt MM: I left out this mysterious remark: "MATLAB was able to find the actual roots when the number of observations are greater than or equal to one million. In the case of hundred thousand and five hundred thousand observations, the function failed to find the actual roots". Don't the experiments show that for these values you did find roots?} 
The figures show that, when increasing the number of observations $m$, the algorithm increasingly well succeeds in identifying the structure of the underlying routing matrix. \michel{The estimation in this experiment (and in all subsequent experiments) was performed using the procedure described in (P2) in Subsection~\ref{ssec:momeq}.}

\medskip

{\it Experiment 3.} 
% Here your 10-nodes example, to demonstrate we can in principle deal with large instances. Explain how you picked the parameters; I recall for instance that the entries of $Q$ were somehow sampled? Of course not needed to give the precise values of all entries in ${\bs\theta}$. Then you can give the largest percentagewise error of all entries of ${\bs\theta}$ (hopefully small). Specify $m$ and the computation time (simulation time is irrelevant). 
% PUTTING IN SCREENSHOTS FOR THE BIG EXAMPLE
In this experiment we demonstrate that the estimation procedure still gives satisfactory output even for networks of a relatively high dimension. In the instance considered we let $\bs{p} = \bs{1}$ be known. The transition matrix $Q$ is a uniformly sampled sub-stochastic matrix with diagonal values equal to $0$.
%The transition matrix $Q$ is the following uniformly sampled sub-stochastic matrix:
% \[{\footnotesize
% Q = 
% \begin{pmatrix*}[r]
% 0 & 0.0608 & 0.1291 & 0.0975 & 0.0124 & 0.4925 & 0.0147 & 0.0906 & 0.0339 & 0.0116 \\
% 0.1726 & 0 & 0.0164 & 0.0141 & 0.1440 & 0.0047 & 0.0272 & 0.1264 & 0.0278 & 0.1110 \\
% 0.1464 & 0.1841 & 0 & 0.3158 & 0.1069 & 0.0073 & 0.0211 & 0.0086 & 0.0672 & 0.0582 \\
% 0.1970 & 0.1819 & 0.2441 & 0 & 0.0543 & 0.1741 & 0.0046 & 0.0166 & 0.0038 & 0.1030 \\
% 0.0232 & 0.1105 & 0.3216 & 0.0397 & 0 & 0.0796 & 0.0470 & 0.0507 & 0.2484 & 0.0785 \\
% 0.1476 & 0.1197 & 0.1102 & 0.0299 & 0.1680 & 0 & 0.0212 & 0.0574 & 0.1247 & 0.0432 \\
% 0.0180 & 0.0038 & 0.0963 & 0.1402 & 0.1720 & 0.1457 & 0 & 0.1050 & 0.0106 & 0.2736 \\
% 0.0103 & 0.0488 & 0.1286 & 0.0013 & 0.0543 & 0.0849 & 0.0135 & 0 & 0.2743 & 0.1553 \\
% 0.1196 & 0.1553 & 0.2565 & 0.1716 & 0.0568 & 0.0515 & 0.0978 & 0.0368 & 0 & 0.0157 \\
% 0.0698 & 0.0067 & 0.0121 & 0.0582 & 0.0443 & 0.0814 & 0.0666 & 0.3983 & 0.1942 & 0 \\
% \end{pmatrix*}}
% \]
The service times follow exponential distributions with rate parameters given by 
$
\bs\mu = (2,3,3,6,3,4,6,3,2,5),
$
and the external arrival rates are given by
$
\bs\lambda = (3,2,4,5,2,6,4,3,2,6).$

\begin{figure}
    \centering
    \includegraphics[width=0.65\linewidth]{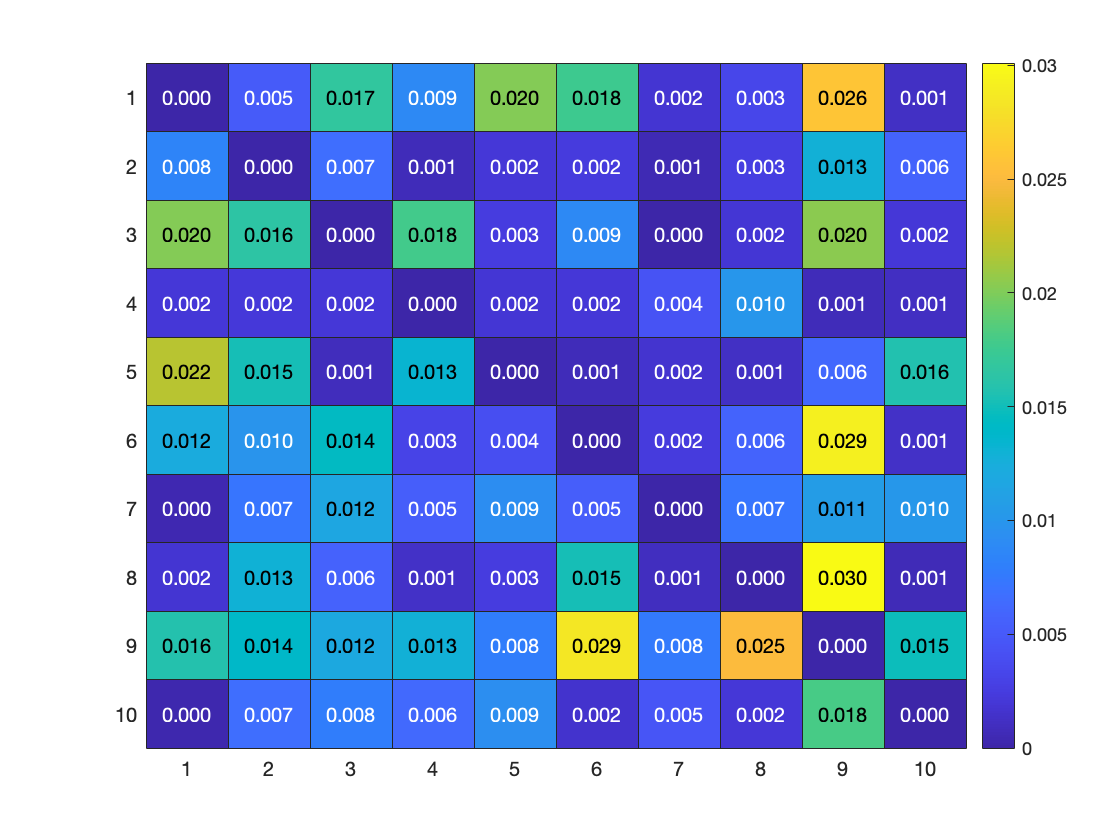}
    \caption{Colormap of $|Q - \hat{Q}|$}
    \label{fig:colormap}
\end{figure}

Using a dataset of $m=2\cdot 10^6$ observations obtained with sampling performed at rate $\beta = 10$, the estimation procedure completed in 25.45 seconds. The estimated values are as follows:
%\footnote{\scriptsize \tt MM: huh?! Above you say service time distributions are known?}
$$ \hat{\bs\mu}  = (1.997,3.015,2.921,5.931,3.019,3.924,5.980,3.000,1.996,5.048) \,.$$
% \[{\footnotesize
% \hat Q = 
% \begin{pmatrix*}[r]
% 0 & 0.0659 & 0.1458 & 0.1064 & 0.0326 & 0.4749 & 0.0165 & 0.0873 & 0.0598 & 0.0128 \\
% 0.1642 & 0 & 0.0231 & 0.0149 & 0.1461 & 0.0071 & 0.0277 & 0.1291 & 0.0407 & 0.1052 \\
% 0.1666 & 0.2005 & 0 & 0.3336 & 0.1094 & 0.0162 & 0.0212 & 0.0105 & 0.0876 & 0.0603 \\
% 0.1950 & 0.1841 & 0.2462 & 0 & 0.0562 & 0.1761 & 0.0002 & 0.0066 & 0.0046 & 0.1019 \\
% 0.0451 & 0.1257 & 0.3207 & 0.0530 & 0 & 0.0804 & 0.0487 & 0.0495 & 0.2546 & 0.0942 \\
% 0.1354 & 0.1295 & 0.1246 & 0.0330 & 0.1721 & 0 & 0.0188 & 0.0516 & 0.1539 & 0.0446 \\
% 0.0184 & 0.0110 & 0.0846 & 0.1455 & 0.1812 & 0.1403 & 0 & 0.0978 & 0.0000 & 0.2838 \\
% 0.0121 & 0.0618 & 0.1229 & 0.0000 & 0.0568 & 0.0697 & 0.0125 & 0 & 0.3044 & 0.1559 \\
% 0.1354 & 0.1692 & 0.2684 & 0.1845 & 0.0646 & 0.0800 & 0.0902 & 0.0619 & 0 & 0.0306 \\
% 0.0701 & 0.0000 & 0.0199 & 0.0645 & 0.0537 & 0.0830 & 0.0712 & 0.3963 & 0.2123 & 0 \\
% \end{pmatrix*}} \,.
% \]
For the estimate of $Q$, the absolute error $|Q - \hat{Q}|$ is shown in Figure~\ref{fig:colormap}. The total squared error across the 90 off-diagonal entries of the transition rate matrix is $|Q - \hat{Q}|^2 =0.0114$.
%Although the zero elements in $Q$ are not estimated exactly as zero, their estimated values remain negligible (with the maximum absolute error being 0.0300).\footnote{\footnotesize \tt MM: I thought you enforced these elements to be equal to 0?!?}
While the estimates of ${\bs{\lambda}}$ show a significant deviation from the true values, the effective arrival rate estimates ${\bs{\lambda}}^{\mathrm{eff}}$ demonstrate much closer agreement. Specifically, we have
\[
\hat{\bs \lambda}^{\rm eff} = 
(
36.089, 34.417, 44.766, 35.1124, 28.759, 40.797, 14.474, 30.051, 33.806, 30.052
) \,,
\]
while
\[
\bs \lambda^{\rm eff} = 
 (36.139, 34.242, 45.978, 35.521, 28.583, 41.585, 14.523, 30.051, 33.867, 29.763
 ).
\]
The above effect is a consequence of the fact that our estimation procedure follows a sequential approach, as described in Subsection \ref{ssec:p1}: we first estimate $\bs \mu$ and $Q$, and subsequently determine $\hat{\bs\lambda}$ based on the estimated parameters $\hat{\bs \mu}$ and $\hat{Q}$. Consequently, the estimation errors potentially propagate. Additional experiments reveal that the estimation accuracy degrades in the number of stations, as anticipated.

\medskip

% For $Q$, some percentage errors are coming out to be $100\%$, because true value is very small and estimated value is $0$. For example, see $Q(10,2)$, true value is $0.0067$ and estimated is $0$. The maximum absolute error in $Q$ is $0.0300$. Maximum error for $\mu$ is $2.64\%$.\\\\

\begin{figure*}[h!]
\centering
\includegraphics[clip, trim=6.5cm 19cm 5.5cm 9cm, width=0.6\textwidth]{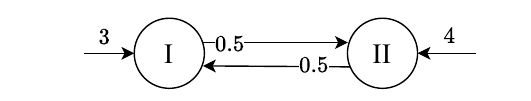}
\caption{Two-station network used to test the model-free approach.}\label{fig:model_free}
\end{figure*}

{\it Experiment 4.}  We now consider an experiment in which we do not know the parametric form of the service-time distribution. The main goal is to show how the model-free approach of Subsection~\ref{ssec:mfree} can be used to produce accurate estimates. We consider an instance with $n=2$ stations, so that in the model-free approach we need to estimate
\[({\mathbb E}[G_1],{\mathbb E}[G_2]),\:\:({\mathscr G}_1(\beta),{\mathscr G}_2(\beta)),\:\:(\dot{\mathscr G}_1(\beta),\dot{\mathscr G}_2(\beta)).\]
In this experiment we take $\bs{p} = \bs{1}$, while the external arrival rates are \( \bs{\lambda} = (3, 4)^\top \) and the transition probabilities $q_{1,2} = q_{2,1} = 0.5, \, q_{1,1} = q_{2,2} = 0$; see Figure~\ref{fig:model_free}.

 We consider two cases: in the former the data is generated such that the service times are exponentially distributed with rates $\bs\mu$, whereas in the latter the service times follow and Erlang-2 distribution with rate $\bs \mu$. We take  \( \bs{\mu} = (3, 5)^\top \), so that the resulting expected service time vectors for the two cases are $(0.3333,\, 0.2)$ and $(0.6667, 0.4)$, respectively.

\begin{table}[h]
\centering
\begin{tabular}{|c|c|c|c|}
  \hline
   & $\hat{\bs \lambda}$ & ${\mathbb E}[G_1],{\mathbb E}[G_2]$ & $\hat{q}_{12}, \, \hat{q}_{21}$ \\
  \hline
  Exponential assumption & $3.0196, \, 3.8594$ & $0.3545, \, 0.2156$ & $0.4695, \, 0.4769$ \\
  \hline
  Model-free & $3.0250, \, 4.0209$ & $0.3330, \, 0.2003$ & $0.4948, \,  0.4977$ \\
  \hline
\end{tabular}
\caption{Service times in data are exponentially distributed with rate $\bs \mu$. }
\label{tab:Exp}
\end{table}
%%%%%%%%%%%%%%%%%%%%%%%%
\begin{table}[h]
\centering
\begin{tabular}{|c|c|c|c|}
  \hline
   & $\hat{\bs \lambda}$ &${\mathbb E}[G_1],{\mathbb E}[G_2]$ & $\hat{q}_{12}, \, \hat{q}_{21}$ \\
  \hline
  Exponential assumption & $2.9695, \, 4.4515$ & $0.5192, \, 0.3050$ & $0.6027, \, 0.5814$ \\
  \hline
  Model-free & $2.8990, \, 3.8911$ & $0.6687, \, 0.4003$ & $0.5162, \, 0.5115$ \\
  \hline
\end{tabular}
\caption{Service times in data are Erlang distributed with parameter $2$ and rate $\bs \mu$.}
\label{tab:Erlang2}
\end{table}
Tables \ref{tab:Exp} and \ref{tab:Erlang2}
% compare estimates in which {\bf the
% sampling are performed at rate $\beta = 5$}, and the service times are assumed to be exponentially distributed, with the model-free approach. 
compare the estimates obtained with sampling performed at rate $\beta = 3$ and service times assumed to be exponentially distributed, with the model-free approach. 
They show that when the underlying service times in the data are exponentially distributed, both approaches produce reliable estimates, with the former achieving slightly higher accuracy. However, when the service times follow an Erlang-2 distribution, the `exponential-based' method suffers from the model misspecification, resulting in estimates that deviate significantly from the true values. In contrast, the model-free approach remains robust and continues to yield accurate results. This highlights a key advantage of the model-free method: its ability to handle uncertainty about the parametric form of service-time distributions effectively.

\section{Concluding remarks}

In this paper we have studied a queueing network in which each node is modeled as an M/G/$\infty$ queue. Based solely on noisy observations of the number of individuals at each node, we have developed a procedure to identify the service rates, arrival rates, and the routing matrix. Notably, we can infer not only the existence of links between nodes but also their direction --- an essential feature for causality analysis.

Our work has practical relevance in multiple domains. For example, in social networks, understanding the routing matrix, and thus the underlying structure, can help identify influencers, communities, or potential spreaders of information (or misinformation). Such insights are valuable in applications such as marketing, healthcare, and political campaigns. Moreover, knowledge of the network structure facilitates the detection of bottlenecks and supports strategic decisions, such as whether to add or remove a route between specific nodes.
%\textcolor{red}{Our framework may also be useful for fault detection. For instance, if a link suddenly fails, this will inevitably affect the estimates of the moments of the node populations. If one could identify the time of the changepoint, one can also identify which link has gone down.}
% Although single-snapshot moments and cross moments can estimate the routing matrix and thereby detect such faults, observing the dynamics of cross moments over time may also reveal structural changes in the network. {\bf Jiesen: Although the routing matrix can be estimated to detect such faults, monitoring the dynamics of empirical cross moments over time may also reveal structural changes in the network.}

Future research directions include:
\begin{itemize}
    \item[$\circ$] {\it Single-server queue analysis}: Extend the analysis to networks where each station is a {\it single}-server queue (operating under, say a first-come, first-served discipline), rather than an infinite-server queue. An even more general setup could concern queues where the expected waiting time increases with the number of individuals already present. In particular, one could consider cases in which the arrival rate at a station can be modeled as a function of its current occupancy, reflecting customer reluctance to join when the queue is more congested. A simple example is a step-function arrival rate: for instance, let (for some \(m_{0,i} \in \mathbb{N}^+\)) the arrival rate \(\lambda_i\) remain as before when \michel{\(M_i \leqslant m_{0,i}\)} and set \(\lambda_i = 0\) otherwise.
    \item[$\circ$] {\it Endogenous routing dynamics}: Study scenarios where customers {\it choose} their next station after service completion, making the routing matrix endogenous rather than exogenous. This could involve modeling customer decisions based on real-time information on queue lengths or other system states.
\end{itemize}

\bibliographystyle{plain}
\bibliography{bibliography.bib}
% \bibliography{Inf_serv_netw.bib}

 \end{document}